\documentclass{amsart}

\usepackage{graphicx}
\usepackage{amssymb}
\usepackage{amscd}
\usepackage{amsmath}
\usepackage{amsfonts}
\usepackage{a4}
\usepackage{amsthm}
\usepackage{color}

\theoremstyle{plain}

\newtheorem{theo}{Theorem}[section] 
\newtheorem{prop}[theo]{Proposition}
\newtheorem{lemme}[theo]{Lemma}
\newtheorem{cor}[theo]{Corollary}
\newtheorem{defin}[theo]{Definition}
\newtheorem{AS}[theo]{Assumption}

\theoremstyle{definition}

\newtheorem{rem}[theo]{Remark}
\newtheorem{exe}[theo]{Example}

\newcommand{\eval}[1]{\lvert_{#1}}

\newcommand{\Ci}{C^{\infty}}
\newcommand{\op}{\operatorname}

\newcommand{\com}[1]{}

\newcommand{\R}{{\mathbb{R}}}
\newcommand{\N}{{\mathbb{N}}}
\newcommand{\C}{{\mathbb{C}}}

\newcommand{\T}{{\mathbb{T}}}

\newcommand{\de}{{\delta}}
\newcommand{\be}{{\beta}}
\newcommand{\al}{{\alpha}}
\newcommand{\la}{{\lambda}}
\newcommand{\De}{{\Delta}}
\newcommand{\si}{{\sigma}}
\newcommand{\hb}{{\hbar}}
\newcommand{\ka}{{\kappa}}

\newcommand{\te}{{\theta}}
\newcommand{\om}{{\omega}}
\newcommand{\Om}{\Omega}
\newcommand{\ph}{{\varphi}}

\newcommand{\Ga}{{\Gamma}}
\newcommand{\eps}{{\epsilon}}

\newcommand{\muquant}{{\op{M}}}

\newcommand{\td}{\mathfrak{t}_d}

\newcommand{\fl}{{\mathcal{F}}_{\la}}
\newcommand{\fr}{{\mathcal{F}}_{r}}
\newcommand{\flr}{{\mathcal{F}}_{\la,r}}

\newcommand{\laT}{{\mathcal{T}}_{\la}}
\newcommand{\rT}{{\mathcal{T}}_{r}}

\newcommand{\p}{{\mathfrak{p}}}

\newcommand{\laH}{{\mathcal{H}}_{\la,k}}

\newcommand{\laP}{\Pi_{\la,k}}
\newcommand{\rP}{\Pi_{r,k}}

\author{L. CHARLES}

\title{Toeplitz operators and Hamiltonian torus actions}   

\address{Institut de Math{\'e}matiques \\ Analyse
  Alg{\'e}brique \\ Universit{\'e} Pierre et Marie Curie \\ 175, rue du
  Chevaleret \\ 75013 Paris \\ FRANCE}

\thanks{This work has been partially supported by the European
Commission under the Research Training Network (Mathematical Aspects of
Quantum Chaos) n° HPRN-CT-2000-00103 of the IHP Programme. }



\begin{document}

\begin{abstract}
This paper is devoted to semi-classical aspects of
symplectic reduction. Consider a compact prequantizable K{\"a}hler
manifold $M$ with a Hamiltonian torus action. In the seminal paper
\cite{GuSt1}, Guillemin and Sternberg introduced an isomorphism between
the invariant part of the quantum space associated to $M$  and the
quantum space associated to the symplectic quotient of $M$, provided
this quotient is non-singular. 
We prove that this isomorphism is a Fourier integral operator and that
the Toeplitz operators of $M$ descend to Toeplitz
operators of the reduced phase space. We also extend these results to
the case where the symplectic quotient is an orbifold and estimate the spectral density of a reduced Toeplitz
operator, a result related to the Riemann-Roch-Kawasaki theorem.  
\end{abstract}

\maketitle

\bibliographystyle{plain}

\section{Introduction}

Consider a symplectic manifold $(M,\om)$ with a Hamiltonian action of
a $d$-dimen\-sional torus $\T^d$. Let $\mu$ be a momentum map. Following
Marsden-Weinstein \cite{MaWe}, if $\la$ is a regular value of $\mu$, the reduced space 
$$ M_r := \mu^{-1}(\la) / \T^d$$
is naturally endowed with a symplectic form $\om_r$. The quantum
analogue of this reduction has been the subject of important studies,
starting from the paper \cite{GuSt1} of Guillemin and Sternberg, and
has led to many versions of the ``quantization commutes with reduction'' theorem. 
In most of these articles, the quantization is defined as a Riemann-Roch number or the index of a spin-c Dirac operator which represents the dimension of a virtual quantum space, cf. the review article \cite{Sj}. The relationships between deformation quantization and symplectic reduction have also been considered \cite{Xu}, \cite{Fe}. 

This paper is devoted to the quantum aspects of symplectic reduction
in the semi-classical setting. Here the quantization consists of a
Hilbert space  with a semi-classical algebra of operators. More
precisely, we assume that $M$ is compact, K{\"a}hler and endowed with a
prequantization bundle $L \rightarrow M$, i.e. a Hermitian
line bundle with a connection of curvature $-i \om$. For every
positive integer $k$, let us define the quantum space
${\mathcal{H}}_k$ as the space of holomorphic sections of $L^k$. The semi-classical limit is $k
\rightarrow \infty$ and the operators we will consider are the Toeplitz
operators, introduced by Berezin in
\cite{Be}. The application of microlocal techniques in this context
started with Boutet de Monvel and Guillemin \cite{BoGu}. This point of
view made it possible to extend many results known for the
pseudodifferential operators with small parameter to the Toeplitz
operators, as for instance, the trace formula \cite{BP2} and the
Bohr-Sommerfeld conditions \cite{oim2}.

Assume that the torus action preserves the complex
structure of $M$. Then following Kostant and Souriau we can associate to the
components $(\mu_1,...,\mu_d)$ of the moment map $\mu$ some commuting
operators $\muquant_1,...,\muquant_d : {\mathcal{H}}_k \rightarrow
{\mathcal{H}}_k$. Suppose that  $\la = (\la_1,...,\la_d)$ is a joint
eigenvalue of these operators when $k=1$ and that the torus action on
$\mu^{-1}(\la)$ is free. The following  theorem is a slight
reformulation of the main result of Guillemin and Sternberg.

\begin{theo}[\cite{GuSt1}] \label{i1}
$M_r$ inherits a natural K{\"a}hler structure and a prequantization bundle
$L_r$, which defines quantum spaces
${\mathcal{H}}_{r,k}$. Furthermore, for any $k$,
there exists a natural vector space isomorphism $\op{V}_k$ from $${\mathcal{H}}_{\la,k} := \cap_{i=1}^{d} \op{Ker} ( \op{M}_i -
\la_i)$$ onto ${\mathcal{H}}_{r,k}$.
\end{theo}

The various quantum spaces have natural scalar products induced by the Hermitian
structure of the prequantization bundles and the Liouville measures,
but unfortunately the isomorphism $ \op{V}_k$ is not
necessarily unitary. So we will use 
$$\op{U}_k :=  \op{V}_k  (
\op{V}_k^* \op{V}_k)^{-\frac{1}{2}}: {\mathcal{H}}_{\la,k} \rightarrow
{\mathcal{H}}_{r,k} $$ instead.
Our first result relates the Toeplitz operators of $M$ with the
Toeplitz operators of $M_r$.

\begin{theo} \label{i2}
Let $(\op{T}_k : {\mathcal{H}}_{k} \rightarrow {\mathcal{H}}_{k}) _{k
  \in \N^*}$ be a Toeplitz operator of $M$ which commutes with $\muquant_1$,
  ...,$\muquant_d$, and with principal symbol $f \in \Ci(M)$. Then 
$$ (\op{U}_k \op{T}_k \op{U}_k^* : {\mathcal{H}}_{r,k} \rightarrow
{\mathcal{H}}_{r,k})_k$$
is a Toeplitz operator of $M_r$. Furthermore, $f$ is $\T^d$-invariant
and the principal symbol
$f_r$ of $ (\op{U}_k \op{T}_k \op{U}_k^*)$ is such that $p^*f_r = j^* f$, where $p$
and $j$ are respectively the projection $\mu^{-1}(\la)
\rightarrow M_r$ and the embedding $\mu^{-1}(\la) \rightarrow M$.
\end{theo} 

When the torus action on $\mu^{-1}(\la)$ is not free but locally free,
the reduced space $M_r$ is not a manifold, but an orbifold. These spaces
with finite quotient singularities were first introduced  by Satake in
\cite{Sa1}. Many results or notions of differential geometry have been
generalized to orbifolds: index theorem \cite{Ka}, fundamental
group \cite{Th}, string theory \cite{Ru}.  
Not surprisingly, theorems \ref{i1} and \ref{i2} are still valid in this case. Motivated by
this, we prove the basic properties of the Toeplitz operators on the
orbifold $M_r$: description of their Schwartz kernel and the symbolic calculus.
Our second main result is the estimate of the spectral density of a
Toeplitz operator on the orbifold $M_r$. This simple result in the
manifold case involves here oscillatory contribution of the
inertia orbifolds or twisted sectors associated to $M_r$ and is related
to the Kawasaki-Riemann-Roch theorem \cite{Ka}, cf. theorems \ref{int2}
and \ref{Tra} for precise statements.

With a view towards application, we also consider the simple case
where $M$ is $\C^n$ with a linear circle action whose momentum map is
a proper harmonic oscillator. The quantum data associated with $\C^n$
are defined by the Bargmann representation and the reduced space is a
twisted projective space. Actually, this is nearly a particular case
of the previous setting, except that $\C^n$ is not compact. As a
corollary of theorem \ref{i2}, the spectral analysis of an operator
commuting with the quantum harmonic oscillator is reduced to that
of a Toeplitz operator on a projective space.  In collaboration with
San Vu Ngoc, we plan to apply this to the semi-excited spectrum of a
Schr{\"o}dinger operator with a non-degenerate potential well.

We also prove that the isomorphism $\op{V}_k$ of theorem \ref{i1} and
its unitarization $\op{U}_k$ are Fourier integral operators. Thus we can interpret theorem \ref{i2} as a composition of Fourier integral
operators with underlying compositions of canonical
relations. Actually our proof of theorem \ref{i2} is elementary
in the sense that it relies on the geometric
properties of the isomorphism $\op{V}_k$ and doesn't use the
usual tools of microlocal analysis. But with the more general point of view of Fourier integral operators, we hope that we can extend the ``quantization commutes
with reduction'' theorems by using microlocal techniques. For instance,
theorems \ref{i1} and \ref{i2} should hold with general Toeplitz
operators $\op{M}_i$ whose joint principal symbol define a momentum map 
without assuming that the action preserves the complex structure. Also
the K{\"a}hler structure is certainly not necessary. This
microlocal approach is also related to another paper of Guillemin and
Sternberg \cite{GuSt3} (cf. sections \ref{laTose} and \ref{laFo} for a
comparison with our results).

The organization of the paper is as follows. Section \ref{stosc} contains
detailed statements of our main results for the harmonic oscillator on
$\C^n$. In section \ref{s2}, we introduce our set-up in the compact
K{\"a}hler case and recall the results of \cite{GuSt1} proving that the
reduced quantum space is isomorphic to the joint eigenspace. Section \ref{rTo}
 contains the statements and proofs of our main results for the
 reduction of Toeplitz operators. In section \ref{vasavoir}, we
 interpret these results as compositions of Fourier integral
 operators. Sections \ref{ToepOrb} and \ref{prep} are devoted to the Toeplitz operators on K{\"a}hler orbifolds.

{\bf Acknowledgment.} I would like to thank Sandro Graffi for his encouragement and Alejandro Uribe for introducing me to the quantization of symplectic reduction a few years ago.


\section{Statement of the results for the harmonic oscillator} \label{stosc}

Assume $\C^n$ is endowed with the usual symplectic 2-form 
$\om = i
(dz_1 \wedge d\bar{z}_1 + ...+ dz_n \wedge d \bar{z}_n).$
Let $H$ be the harmonic oscillator
$$ H := \p_1 |z_1|^2+...+\p_n |z_n|^2 $$ 
where $\p_1,...,\p_n$ are positive relatively prime integers. Consider the scalar product
\begin{gather} \label{defpsbargman}
 (\Psi ,\Psi' )_{\C^n} = \int_{\C^n} e^{- \hb^{-1} |z|^2} \Psi (z)
 . \bar{\Psi}'(z) \tfrac{|dz.d\bar{z}|}{n!} 
\end{gather} 
where $\Psi,\Psi'$ are functions on $\C^n$ and $|z|^2 = |z_1|^2+...+
|z_n|^2$. The Bargmann space ${\mathcal{H}}$ is the Hilbert space of
holomorphic functions $\Psi$ on $\C^n$ such that  $ (\Psi ,\Psi )_{\C^n} < \infty$. The quantum harmonic oscillator is the unbounded operator of ${\mathcal{H}}$
\begin{gather} \label{defHquant}
 \op{H} :=   \hbar( \p_1 z_1 \partial_{z_1}+ ...+ \p_n z_n \partial_{z_n} ) \end{gather}
with domain the space of polynomials on $\C^n$. 

\subsection{Symplectic reduction} \label{credo}

The Hamiltonian flow of $H$ induces an action of $S^1$ on the level set $P:= \{ H =1 \}$ 
$$  S^1 \times P \rightarrow P, \quad \te,z \rightarrow
l_\te.z =(z_1 e^{i\te \p_1},...,z_n e^{i\te \p_n}) \text{ if } z = (z_1,...,z_n). $$
Define the reduced space $M_r$ as the quotient $P/S^1$. If $\p_1 =
...=\p_n=1$, the action is free, $M_r$ is a manifold and the
projection $P \rightarrow M_r$ is the Hopf fibration. When the $\p_i$
are not all equal to $1$, the action is not free, but locally
free. Hence $M_r$ is not a manifold, but an orbifold. In any cases,  $M_r$ is naturally endowed with a symplectic $2$-form $\om_r$.

We may also define a complex structure on the space $M_r$ by viewing it as a complex quotient. Consider the holomorphic action of $\C^*$ 
\begin{gather} \label{accoCn}
 \C^* \times \C^n \rightarrow \C^n, \quad u,( z_1,...,z_n) \rightarrow
 (z_1 u^{\p_1},...,z_n u^{\p_n} ) \end{gather} 
Each $\C^*$-orbit
 of $\C^n -\{ 0\}$ intersects $P$ in a $S^1$-orbit, which identifies
 $M_r$ with $\C^n -\{ 0\} / \C^*$. This quotient is called a twisted
 projective space, the standard projective space is obtained when $\p_1 = ...=\p_n=1$. The complex structure is compatible with the symplectic form $\om_r$. So $M_r$ is a K{\"a}hler orbifold.

\subsection{Quantum Reduction}
$\op{H}$ has a discrete spectrum given by
$$ \op{Sp}(\op{H}) = \{  \hbar ({\mathfrak{p}}_1 \al(1) + ...+
{\mathfrak{p}}_n \al(n)) ; \; \al \in \N^n \}$$
Hence $1$ is an eigenvalue only if $\hb$ is of the form $1/k$ with $k$
a positive integer. Since we will only consider the eigenvalue $1$, we
assume from now on that $$\hb = 1/k  \text{ with } k \in \N^*$$ and
use the large parameter $k$ instead of the small parameter $\hb$. We
denote by ${\mathcal{H}}_k$ and $\op{H}_k$ the Bargmann space and the
quantum harmonic oscillator. The vector space $${\mathcal{H}}_{1,k} :
= \op{Ker}( \op{H}_k -1)$$ is generated by the monomials $z^{\al}$
such that ${\mathfrak{p}}_1 \al(1) + ...+ {\mathfrak{p}}_n \al(n) =
k$. So a state $\Psi \in {\mathcal{H}}_k$ belongs to
${\mathcal{H}}_{1,k}$ if and only if it is invariant in the sense
that  
$$ \Psi ( z_1 u^{\p_1},...,z_n u^{\p_n} ) = u^k \Psi( z_1 ,...,z_n) .$$
Hence there is a holomorphic line orbi-bundle $L_r \rightarrow
M_r$ such that $\mathcal{H}_{1,k}$ identifies with the space
${\mathcal{H}}_{r,k}$ of holomorphic sections of $L_r^k$. $L_r$ has a
natural Hermitian structure and connection of curvature $-i\om_r$,
which turns it into a prequantization orbi-bundle (cf. section \ref{s2}). We
denote by $\op{V}_k $ the isomorphism from ${\mathcal{H}}_{1,k}$ to ${\mathcal{H}}_{r,k}$

\subsection{Reduction of the operators} \label{redosc}

On the Bargmann space a usual way to define operators is the Wick or Toeplitz quantization. Denote by $\Pi_k$ the orthogonal projector of $L^2( \C^n, e^{-k |z|^2}  |dz.d\bar{z}| )$ onto ${\mathcal{H}}_k$. To every function $f$ of $\C^n$ we associate the operator $\op{Op} (f)$ of ${\mathcal{H}}_k$ defined by
$$ \op{Op} (f): \Psi \rightarrow \Pi_k( f.\Psi )$$
More generally, we consider multiplicators $f$ which depend on $k$. Define the class $S(\C^n)$ of symbols $f(.,k)$ which are sequences of $\Ci(\C^n)$ satisfying
\begin{itemize} 
\item there exists $C >0$ and $N$ such that $| f(z, k) | \leqslant C ( 1 + |z| )^N, \; \forall z \in \C^n, \; \forall k$. 
 \item  $f(., k)$ admits an asymptotic expansion of the form $$\textstyle{\sum}_{l=0}^{\infty}  k^{-l} f_l + O(k^{-\infty})$$ with $f_0,f_1,..\in \Ci(\C^n)$ for the $\Ci$ topology on a neighborhood of $P$.
\end{itemize}
For such a symbol, we consider $\op{Op}(f(.,k))$ as an unbounded
operator of ${\mathcal{H}}_k$ with domain polynomials on
$\C^n$.  Its {\em principal symbol} is the function $f_0$. If $f(.,k)$ is invariant with respect to the Hamiltonian flow of $H$, then $\op{Op}(f(.,k))$ sends ${\mathcal{H}}_{1,k}$ into itself.
\begin{rem} \label{expol}
The class of Toeplitz operators with symbol in $S(\C^n)$ contains the
algebra of differential operators generated by $ \tfrac{1}{k}
\partial_{z_i}$  and $z_i$. Indeed, let 
$$ f(.,k) = P_0 + k^{-1}P_1 + ...+ k^{-M} P_M$$ 
where $P_0(\bar{z},z),..,P_M(\bar{z},z)$ are polynomials of $\C [\bar{z},z ]$. Then $\op{Op}(f(.,k))$ is the operator 
$$   P_0 ( \tfrac{1}{k} \partial_z, z) + k^{-1}  P_1 (\tfrac{1}{k}
\partial_z, z)+...+ k^{-M}  P_M ( \tfrac{1}{k} \partial_z ,z). $$
Its principal symbol is $P_0$. If the $P_i$ are linear
combinations of the monomials $\bar{z}^\al z^\be$ such that $\langle
{\mathfrak{p}}, \al - \be \rangle =0$, 
then $f(.,k)$ Poisson commutes with $H$ and $\op{Op}(f(.,k))$
preserves the eigenspaces of the quantum harmonic oscillator.
\qed
\end{rem}

Let $\op{U}_k$ be the unitary map $\op{V}_k  (
\op{V}_k^* \op{V}_k)^{-\frac{1}{2}}: {\mathcal{H}}_{1,k} \rightarrow
{\mathcal{H}}_{r,k}$, that we extend to ${\mathcal{H}}_k$ in such a
way that it vanishes on the orthogonal space to
${\mathcal{H}}_{1,k}$. 
The {\em reduced} operator of $\op{Op}(f(.,k))$ is the operator 
$$  \op{U}_k \op{Op}(f(.,k)) \op{U}_k^*: {\mathcal{H}}_{r,k}
\rightarrow  {\mathcal{H}}_{r,k} $$
Our main result says it is a Toeplitz operator. 
 
\begin{theo} \label{redToeposc}
Let $f(.,k)$ be a symbol of $S(\C^n)$. Then there exists a sequence $g(.,k)$ of $\Ci(M_r)$, which admits an asymptotic expansion of the form $\sum_{l=0}^{\infty} k^{-l} g_l + O(k^{-\infty})$ for the $\Ci$ topology, such that 
$$  \op{U}_k\op{Op}(f(.,k)) \op{U}_k^*  = \op{\Pi}_{r,k} g(.,k) + O(k^{-\infty}) $$ 
where $\op{\Pi}_{r,k}$ is the orthogonal projector onto ${\mathcal{H}}_{r,k}$ and the $ O(k^{-\infty})$ is for the uniform norm. Furthermore, the principal symbol $g_0$ of the reduced operator is given by 
$$ g_0(p(x)) = \int_{S^1} f_0(l_{\te}. x) \ \tfrac{|d\te|}{2 \pi}, \quad \forall x \in P$$
where $p$ is the projection $P \rightarrow M_r$ and $f_0$ is the principal symbol of $\op{Op}(f(.,k))$.
\end{theo}

\subsection{Spectral density} \label{spd}
Consider a self-adjoint Toeplitz operator $(\op{T}_k)_k$ of $M_r$,
\begin{gather*} 
  \op{T}_k = \op{\Pi}_{r,k} g(.,k) + O(k^{-\infty}): {\mathcal{H}}_{r,k} \rightarrow  {\mathcal{H}}_{r,k}
\end{gather*}
where $g(.,k)$ is a  a sequence of $\Ci(M_r, \R)$ with an asymptotic expansion of the form  $\sum_{l=0}^{\infty} k^{-l} g_l + O(k^{-\infty})$ in the $\Ci$ topology. Let $d_k$ be the dimension of  ${\mathcal{H}}_{r,k}$ and 
$$\la_{1}(k) \leqslant \la_2(k) \leqslant ... \leqslant \la_{d_k}(k)$$ 
be the  eigenvalues of $\op{T}_k$ counted with multiplicity. 

Let $f \in \Ci(\R)$. The
estimate of $\sum_{i=1}^{d_k} f \bigl( \la_i (k) \bigr)$ as $k \rightarrow
\infty$ is a standard semi-classical result when $M_r$ is a
manifold. In the orbifold case, this result involves the singular
locus of $M_r$. Denote by $G$ the set of $\zeta =e^{i\te} $
such that $$ P_\zeta : = \{ z \in P; \; l_\te.z = z \}$$ is not empty. A
straightforward computation leads to 
\begin{gather*} 
 G = \{ \zeta \in \C^*; \; \zeta^{\mathfrak{p}_i} = 1 \text{ for some
}i \} 
\end{gather*} 
and 
$$P_\zeta = \C_{\zeta} \cap P \text{ with }\C_{\zeta} = \{ z \in \C^n
 ; \; z_i =0 \text{ if } \zeta ^{p_i} \neq 1 \}.$$ 
The Hamiltonian flow of $H$ preserves $P_\zeta$. Let $M_\zeta$ be the
 quotient of $P_\zeta$ by the induced $S^1$-action. It is a twisted projective
 space which embeds into $M_r$ as a symplectic suborbifold. Denote by
 $n(\zeta)$ its complex dimension. Finally, let $m(\zeta)$ be the
 greatest common divisor of $\{ {\mathfrak{p}}_i; \; \zeta^{{\mathfrak{p}}_i}=1\}$.

\begin{theo} \label{int2}
 For every function $f\in \Ci( \R)$, 
\begin{gather*} 
 \sum_{i=1}^{d_k} f \bigl( \la_i (k) \bigr) =     \sum_{ \zeta \in G}
 \Bigl( \frac{k}{2\pi} \Bigr)^{n(\zeta )}  \zeta^{-k}
 \sum_{l=0}^{\infty}  k^{-l} I_l( \zeta ) + O(k^{-\infty})
\end{gather*} 
The leading coefficients are given by 
$$ I_0( \zeta ) = \frac{1}{m(\zeta)} \Biggl( \prod_{i; \;
  \zeta^{{\mathfrak{p}}_i} \neq 1} (1 - \zeta^{p_i} )^{-1} \Biggr)  \int_ {M_{\zeta}}  f(g_0)\  \delta_{M_\zeta}, $$ 
where $ \delta_ {M_\zeta}$ is the Liouville measure of $M_{\zeta}$.
\end{theo}

Observe that $M_1 =M_r$. The other $M_\zeta$ are of positive
codimension and are the closures of the singular stratas of $M_r$. 
Hence at first order, the formula is the same as in the manifold case
$$ \sum_{i=1}^{d_k} f(\la_i) =  \Bigl( \frac{k}{2\pi} \Bigr)^{n-1}
\int_ {M_r} f(g_0) \delta_{M_r} + O(k^{n-2}) $$
Furthermore applying this result with $f \equiv 1$, we obtain an estimate of the
dimension of ${\mathcal{H}}_k$. When $k$ is sufficiently large, this
dimension is also given by the Riemann-Roch-Kawasaki theorem and both
results are in agreement (cf. remark \ref{remRKK}).


\section{The Guillemin-Sternberg isomorphism} \label{s2}

Let $M$ be a compact connected K{\"a}hler manifold. Denote by $\om \in
\Om^{2}(M,\R)$ the fundamental two-form. Assume that $M$ is endowed
with a prequantization bundle $L \rightarrow M$, that is $L$ is a
Hermitian line bundle with a connection of curvature $-i\om$. 
$(M,\om)$ is a symplectic manifold and represents the classical phase
space. For every positive integer $k$ define the quantum space
${\mathcal{H}}_k$ as the space of holomorphic sections of $L^k
\rightarrow M$. 

Assume that $M$ is endowed with an effective  Hamiltonian torus action
\begin{gather} \label{acM} 
 \T^d \times M \rightarrow M, \qquad \te, x \rightarrow l_{\te} .x
\end{gather} 
which preserves the complex structure. Let $\td$ be the Lie algebra of $\T^d$. If $\xi \in \td$, we denote by $\xi^\#$ the associated vector field of $M$. Let $$\mu : M \rightarrow \td^*$$ be the moment map, so
$\om(\xi^\#,.) + d \langle \mu, \xi \rangle = 0$.  

Following Kostant and Souriau, for every $\xi \in \td$  and every
positive integer $k$ we define the operator $\muquant_{\xi,k}$:  
$$ \muquant_{\xi,k} :=  \langle \mu , \xi \rangle +\tfrac{1}{ik}
\nabla_{\xi^{\#}} : {\mathcal{H}}_k \rightarrow {\mathcal{H}}_k. $$
It has to be considered as the
quantization of the classical observable $\langle \mu, \xi \rangle \in
\Ci(M)$. Since the Poisson bracket of $\langle \mu, \xi \rangle$ and
$\langle \mu, \xi' \rangle$ vanishes, one proves that
$\muquant_{\xi,k}$  commutes with $\muquant_{\xi',k}$.     The {\em joint spectrum} of the $\muquant_{\xi,k}$ is the set of covectors $\la \in \td^*$ such that 
$$ \mathcal{H}_{\la,k} := \{  \Psi \in {\mathcal{H}}_k; \quad
\muquant_{\xi,k} \Psi = \langle \la , \xi \rangle \Psi, \; \forall \
\xi \in \td \} $$
is not reduced to $(0)$. 

The joint eigenvalues are related to the values of $\mu$ in the
following way. First, recall the convexity theorem of Atiyah
\cite{At1} and Guillemin-Sternberg \cite{GuSt2}: the image under $\mu$
of the fixed point set of $M$ is a finite set $$ \{ \nu_1,...,\nu_s \}
$$ and $\mu(M)$ is the convex hull of this set. 

\begin{theo}  \label{neccond}
Let $\nu$ be a value of $\mu$ at some fixed point. Let $(\la, k) \in
\td^* \times \N^*$. Then $\la$ belongs to the joint spectrum of the
$\muquant_{\xi,k}$ only if 
\begin{gather} \label{BS}
 \la \in  \mu(M) \cap \bigl( \nu + \tfrac{2\pi}{k} K \bigr) \end{gather}
where $K$ is the integer lattice of $\td^*$.
\end{theo}
The condition (\ref{BS}) doesn't depend on the choice of
$\nu$: since $(M, \om)$ is
endowed with a prequantization bundle, it is known that for every $i,j$
\begin{gather} \label{grgr}
\nu_i - \nu_j \in  2\pi K. \end{gather}
That $ \la \in  \mu(M)$ is necessary has been proved by Guillemin and Sternberg (cf. theorem
5.3 of \cite{GuSt1}). The second condition, $\la \in \nu +
2\pi k^{-1} K$, is an exact Bohr-Som\-merfeld con\-dition, which
follows from the theory of Kostant and Souriau. 

\begin{exe}
Let $M$ be the projective space $\C {\mathbb{P}} ^3$ with $\om$ the Fubiny-Study form
$$ \om = - i \partial \bar{\partial} ( |z_1|^2 + ...+  |z_4|^2), \quad [z_1,...,z_4] \in  \C {\mathbb{P}} ^3 $$
and $L$ the tautological bundle. Consider the torus action
$$ (\te_1, \te_2 ), [z_1,...,z_4] \rightarrow [ z_1 e^{-2i\pi(\te_1 +\te_2)}, z_2 e^{-6i\pi \te_1}, z_3 e^{-6i \pi \te_2}, z_4] $$
with momentum map $\mu = \frac{2\pi}{|z|^2} ( |z_1|^2 + 3 |z_2|^2, |z_1|^2 + 3 |z_3|^2 )$. 
\begin{figure}[h]
\includegraphics[width=10.cm,angle=0]{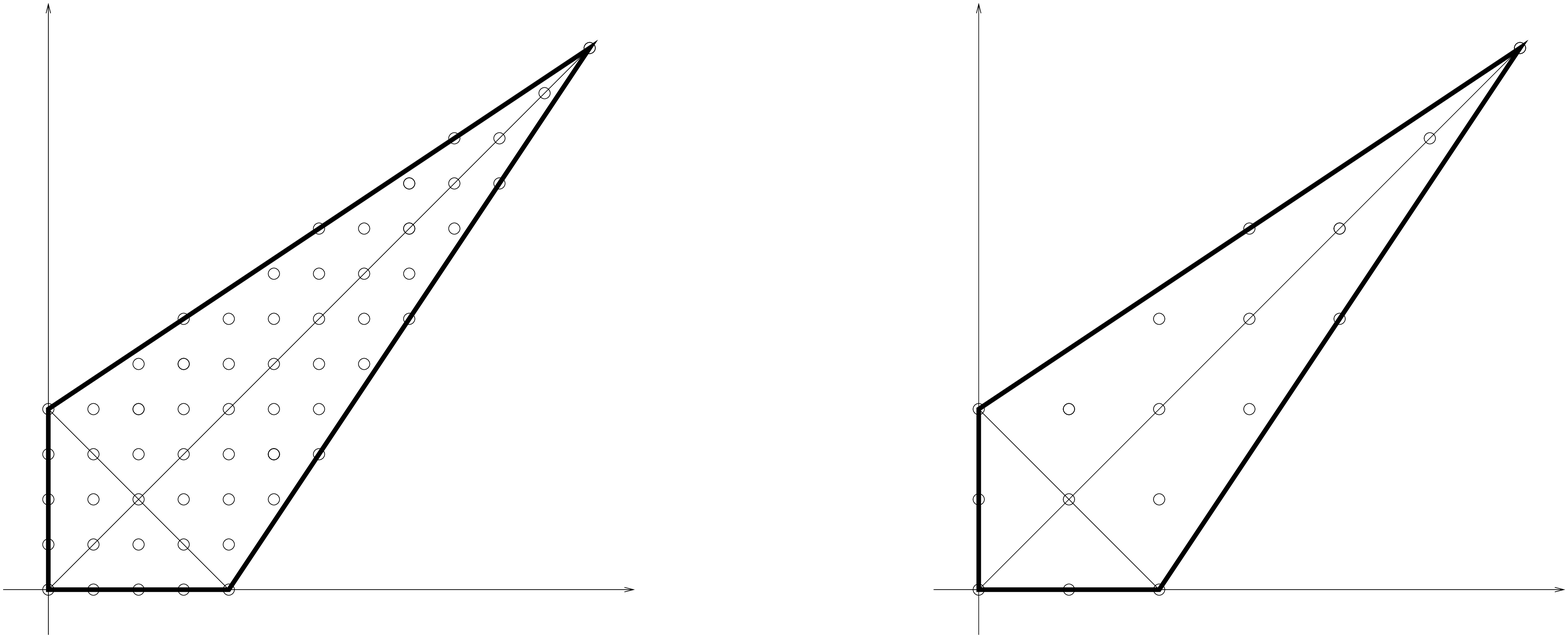}
\end{figure}
The points on the figure are the $\la$ satisfying condition (\ref{BS})
with $k=4$ on the left and $k=2$ on the right. The lines are the critical values of $\mu$. 
\qed \end{exe}

Let $(\la, k) \in \td^* \times \N ^*$. Assume that $(\la, k)$ satisfies
condition (\ref{BS}) and that $\la$ is a regular value of
$\mu$. Denote by $P$ the level set $\mu^{-1}(\la)$. It is known that
$P$ is connected (\cite{GuSt2}, \cite{At1}). The torus
action restricts to a locally free action on $P$. So the quotient
$M_r$ of $P$ is a compact connected orbifold. It is
naturally endowed with a symplectic form $\om_r$.

\begin{theo}[Guillemin-Sternberg \cite{GuSt1}] \label{GuStTh}
$M_r$ inherits by reduction a K{\"a}hler structure with fundamental 2-form $\om_r$ and a prequantization orbi-bundle $L^k_r \rightarrow M_r$ with curvature $-i k \om_r$. Furthermore, there exists a natural isomorphism of vector space 
$$\op{V}_k : {\mathcal{H}}_{\la,k} \rightarrow {\mathcal{H}}_{r,k}$$
where ${\mathcal{H}}_{r,k}$ is the space of holomorphic sections of
$L^k_r$. 
\end{theo}

Now consider a fixed regular value $\la$ of $\mu$ such that the set of
integers $k$ satisfying condition (\ref{BS}) is not empty. Then there
exists a positive integer $\kappa$ such that $(\la, k)$ satisfies
(\ref{BS}) if and only if $k$ is a positive multiple of
$\ka$. Furthermore the K{\"a}hler structure of $M_r$ doesn't depend on
$k$ and for every such $k$, 
$$L^k_r = (L^\ka_r)^{\otimes k/\ka}.$$ 
If $M_r$ is a manifold, it follows from Kodaira vanishing
theorem and Riemann-Roch theorem that 
\begin{gather} \label{estRRK}
\op{dim} {\mathcal{H}}_{r,k} = \Bigl( \frac{k}{2 \pi} \Bigr) ^{n_r} \op{Vol}
(M_r) + O(k^{n_r -1})
\end{gather} 
as $k$ goes to infinity, where $  \op{Vol}
(M_r)$ is the symplectic volume of $M_r$ and $n_r$ its dimension. This gives a partial converse
to theorem $\ref{neccond}$: if $k$ is a sufficiently large positive multiple of
$\kappa$, the eigenspace ${\mathcal{H}}_{\la,k}$ is not reduced to
$(0)$, so  $\la$ belongs to the joint spectrum of the $M_{\xi,k}$. 

If $M_r$ is an orbifold, the same result holds and follows from
Riemann-Roch-Kawasaki theorem \cite{Ka}. Indeed, we assumed that the
torus action is effective. This implies that its restriction to $P$ is also
effective. So $M_r$ is a reduced orbifold or equivalently its
principal stratum has multiplicity one. This explains why formula
(\ref{estRRK}) remains unchanged, without a sum of oscillatory terms.

All the semi-classical results we prove in this paper are in this
regime,
$$ \text{ \em{$\la$ is fixed and $k \rightarrow \infty $ running
    through the set of positive multiples of
$\ka$}} .$$

In the remainder of this section, we recall the main steps of the proof of the
Guillemin-Sternberg theorem. We follow the presentation given by
Duistermaat in \cite{Du2}, that is we consider separately the reduction
of the symplectic and prequantum data and the reduction of the
complex structure. We explain in remarks how the same constructions
apply to the harmonic oscillator.

\begin{rem} (harmonic oscillator).
The Bargmann space can be viewed as a space of holomorphic sections of a prequantization bundle over $\C^n$. Let $L := \C^n \times \C$ be the trivial bundle over $\C^n$. We identify the sections of $L^k$ with the functions on $\C^n$. Introduce a connection and a Hermitian structure on $L^k$ by setting 
$$ \nabla \Psi =  d \Psi - k\Psi (\bar{z}_1 dz_1 +...+\bar{z}_n dz_n) , \qquad (\Psi,\Psi)(z) = e^{- k|z|^2} |\Psi(z)|^2.$$
In this way $L^k$ becomes a prequantization bundle with curvature $-ik
\om$. The scalar product defined in (\ref{defpsbargman}) is  
$$ (\Psi, \Psi)_{\C^n} = \int_{\C^n}  (\Psi,\Psi)(z)  \tfrac{|dz.d\bar{z}|}{n!} $$ 
So the Bargmann space ${\mathcal{H}}_k$ is the Hilbert space of
holomorphic sections $\Psi$ of $L^k$ such that $ (\Psi, \Psi)_{\C^n}$
is finite.  Furthermore it is easily checked that the quantum harmonic
oscillator $\op{H}_k$ defined in (\ref{defHquant}) is given by 
$$ \op{H}_k \Psi = ( H  +\tfrac{1}{ik}  \nabla_{X_H}) \Psi $$
where $X_H$ is the Hamiltonian vector field of $H$. \qed
\end{rem}

\subsection{Reduction of the symplectic and prequantum data} 

We first lift the torus action (\ref{acM}). If $z\in L^k_x$ and $\xi \in
\td$, we denote by ${\mathcal{T}}_{\xi} .z$ the parallel transport of
$z$ along the path 
$$[0,1] \rightarrow M, \quad   s \rightarrow l_{\exp ( s \xi)}.x $$ If $\nu$ is the value of $\mu$ at some fixed point and $\xi$ belong to
the integer lattice of $\td$,
$$ e^{ik \langle \nu - \mu  , \xi \rangle } {\mathcal{T}}_{\xi} .z =
z, \text{ for every $z \in L^k$.} $$
Indeed, this is obviously true if $z \in L^k_x$ where $x$ is a fixed
point and $\mu(x) = \nu$. By proposition 15.3 of \cite{Du1}, the result follows for every $z$.

Consider now $(\la, k) \in \td^* \times \N ^*$ which satisfies
condition (\ref{BS}). Then the action of $\T^d$ on $M$ lifts to $L^k$
$$ \T^d \times L^k \rightarrow L^k, \quad (\te , z) \rightarrow
{\mathcal{L}}_\te .z := e^{ik \langle \la - \mu  , \xi \rangle }
{\mathcal{T}}_{\xi} .z $$
with $\xi \in \td$ such that $\exp \xi = \te$.
One can check the following facts: for
every $\te$, ${\mathcal{L}}_\te$ is an automorphism of the
prequantization bundle $L^k$, it preserves the complex
structure. Furthermore, the obtained representation of $\T^d$ on ${\mathcal{H}}_k$ induces the representation of the Lie algebra $\td$ given by the operators
$$ \nabla_{\xi^{\#}} + i k \langle \mu - \la , \xi \rangle = ik (\muquant_{\xi,k} -  \langle \la , \xi \rangle), \quad \xi \in \td. $$ 
Hence the joint eigenspace ${\mathcal{H}}_{\la,k}$ is the space of invariant holomorphic sections
$${\mathcal{H}}_{\la,k} = \bigl\{ \Psi \in {\mathcal{H}}_k;\;
{\mathcal{L}}^*_\te \Psi =\Psi, \; \forall \ \te \in \T^d \bigr\} .$$

Denote by $j : P \rightarrow M$ and $p : P \rightarrow M_r$ the
natural embedding and projection. Recall that the reduced symplectic
$2$-form $\om_r$ is defined by $ p^*\om_r = j^* \om$.  
Let $L^k_r$ be the quotient of $j^*L^k$ by the torus action. This is an Hermitian orbi-bundle over $M_r$ and $p^*L^k_r$ is naturally isomorphic with $j^*L^k$. Furthermore $L^k_r$ admits a connection $\nabla $ such that $$p^* \nabla = j^* \nabla.$$ Its curvature is $-i k\om_r$. So $L^k_r$ is a prequantization orbi-bundle. 
Since the sections of ${\mathcal{H}}_{\la,k}$ are invariant, their restrictions to $P$ descend to $M_r$, 
\begin{gather} \label{def2}
{\mathcal{H}}_{\la,k} \rightarrow \Ci(M_r ,L^k_r), \quad \Psi \rightarrow \Psi_r \text{ such that } p^*\Psi_r = j^* \Psi .\end{gather}
This is the first definition of the Guillemin-Sternberg
isomorphism. In the case the action is not free and $M_r$ is an
orbifold, more details will be given in remark \ref{orbidata}.

\begin{rem} (harmonic oscillator). As in section \ref{stosc}, we only
  consider the eigenvalue $\la =1$. The lift of the $S^1$-action is
  explicitly given by 
$$ S^1 \times (\C^n \times \C) \rightarrow \C^n \times \C, \quad \te,
(z_1,...,z_n,v) \rightarrow (e^{i\p_1 \te}z_1,..., e^{i\p_n
  \te}z_n,e^{ik\te} v). $$
${\mathcal{H}}_{1,k}$ consists of the invariant holomorphic sections of $L^k$.
If $\Psi$ is such a section, one can check by direct
computations that $(\Psi,\Psi)$ and $\nabla \Psi$ are invariant and
that $\nabla_{X_H} \Psi$ vanishes over the level set $P := \{ H =1
\}$. So the quotient $L^k_r$ of $L^k \eval{P}$ by $S^1$ inherits a
structure of prequantization bundle with curvature $-i k \om_r$. \qed
\end{rem}

\subsection{Complex reduction} \label{comcomp}

Let $\T^d_\C$ be the complex Lie group $\T^d \oplus i \td $ with Lie algebra $\td \oplus i \td$. We consider $\T^d$ as a subgroup of $\T^d_{\C}$. Since the torus action preserves the complex structure of $M$, it can be extended in a unique way to a holomorphic action of $\T^d_\C$
$$(\T^d \oplus i\td) \times M \rightarrow M, \quad (\te +it), x \rightarrow l_{\te +it}.x $$
To do this, for every $\xi \in \td$, we define the infinitesimal
generator of $i \xi$ as $J \xi^{\#}$, where $J$ is the complex
structure of $M$. Since $M$ is compact, we can integrate $J
\xi^{\#}$. Then one can check that this defines a holomorphic action.
 
In a similar way, the action on $L^k$ extends to a holomorphic action 
$$ (\T^d \oplus i\td) \times L^k \rightarrow L^k, \quad (\te +it), z \rightarrow {\mathcal{L}}_{\te +it}.z $$
where ${\mathcal{L}}_{\te+it}$ is an automorphism of complex bundle
which lifts $l_{\te +it}$. This gives a representation of $\T^d_\C$ on
${\mathcal{H}}_k$. The induced representation of the Lie algebra $\td \oplus i \td$ is given by the operators
\begin{gather} \label{repLiecomp}
 \nabla_{\xi^\# + J \eta^{\#}} + i k \langle \mu - \la, \xi + i \eta \rangle, \qquad \xi + i \eta \in \td \oplus i \td .\end{gather}
Furthermore ${\mathcal{H}}_{\la,k}$ is the space of $\T^d_\C$-invariant holomorphic sections.

\begin{rem} (harmonic oscillator). 
The holomorphic action of $\C^*$ was given in (\ref{accoCn}). It lifts
to 
$$ \C^* \times (\C^n \times \C) \rightarrow \C^n \times \C, \quad u,
(z_1,...,z_n,v) \rightarrow (u^{\p_1}z_1,..., u^{\p_n}z_n,u^k v). 
$$
The invariant sections of ${\mathcal{H}}_{k}$ are obviously the sections of
${\mathcal{H}}_{1,k}$ \qed
\end{rem}

Let $P_\C$ be the saturated set $\T^d_\C.P$ of $P$. It is an open set
of $M$. The next step is to consider the quotient of $P_\C$ by
$\T^d_\C$. This have to be done carefully because the $\T^d_\C$-action
on $P_\C$ is not proper. Actually, the map 
\begin{gather} \label{diffeo}
  \td \times P \rightarrow P_\C, \quad t,y  \rightarrow l_{it}.y 
\end{gather}  
is a diffeomorphism. So every $\T^d_\C$-orbit of $P_\C$ intersects $P$
in a $\T^d$-orbit and the injection $P \rightarrow P_\C$ induces a
bijection from $M_r$ onto $P_\C /  \T^d_{\C}$. Furthermore every slice
$U \subset P$ for the $\T^d$-action on $P$ is a slice for the
$\T^d_\C$-action. Viewed as a quo\-tient by a holomorphic
action, the orbifold $M_r$ inherits a complex structure. This complex structure is compatible with $\om_r$.

Similarly, the bundle $L^k_r$ may be considered as the quotient
of $L^k \eval{P_\C}$ by the complex action and inherits a
holomorphic structure. This is the unique holomorphic structure
compatible with the connection and the Hermitian product. Denote by
$p_\C$ the projection $P_\C \rightarrow M_r$ and observe that
$p_\C^{*} L^k_r$ is naturally isomorphic with $L^k \eval{P_\C}$.

The interest of viewing $L^k_r$ and $M_r$ as complex quotients is that
there is a natural identification of the $\T^d_\C$-invariant
holomorphic sections $\Psi$ of $L^k \rightarrow P_\C$ with the
holomorphic sections $\Psi_r$ of $L^k_r \rightarrow
M_r$, given by $p_\C^*\Psi_r = \Psi$. So the map (\ref{def2}) takes
its values in the space ${\mathcal{H}}_{r,k}$ of holomorphic sections
of $L^k_r$.

\begin{defin} \label{defGSmap}
$\op{V}_k : {\mathcal{H}}_{\la,k} \rightarrow {\mathcal{H}}_{r,k}$ is the map which sends $\Psi$ into the section $\Psi_r$ such that  $$ p_\C^* \Psi_r = \Bigl( \frac{2\pi}{k} \Bigr)^{\frac{d}{4}} \Psi \eval{P_\C}\quad \text{ or equivalently } \quad  p^* \Psi_r = \Bigl( \frac{2\pi}{k} \Bigr)^{\frac{d}{4}} j^* \Psi.$$ 
\end{defin}

The rescaling by $(2\pi/k)^\frac{d}{4}$ is such that $\op{V}_k$ and its inverse are bounded independently of $k$ (cf. proposition \ref{rough}).

\begin{rem} \label{orbidata} (Orbifold). Let us detail the previous constructions when the $\T^d$-action is
  not free.  For the basic definitions of the theory of orbifolds,
our references are the section 14.1 of \cite{Du1} and the appendix of
\cite{ChRu}.

As topological spaces, $M_r$ and $L^k_r$ are the quotients $P /
\T^d$ and $j^*L^k / \T^d$. $M_r$ is naturally endowed with a
collection of orbifold charts in the following way. 
 Let $x \in P$, $G
\subset \T^d$ be its isotropy subgroup and $U\subset P$ be a slice at
$x$ for the $\T^d$-action. Denote by $\pi_U$ the projection $U \rightarrow M_r$ and by $|U|
\subset M_r$ its image. Then $(|U|,U,G,\pi_U)$ is an orbifold chart
of $M_r$, i.e. $|U|$ is an open set of $M_r$, $U$ a manifold, $G$ a
finite group which acts on $U$ by diffeomorphisms and $\pi_U$ factors
through a homeomorphism $U/G
\rightarrow |U|$.
These charts cover $M_r$ and satisfy some
compatibility conditions, which defines the orbifold structure of
$M_r$. 

For every such chart, the bundle $L^k$
restricts to a $G$-bundle $$L^k_{r,U} \rightarrow U.$$ 
These bundles are orbifold charts of the orbi-bundle $ L^k_r \rightarrow M_r$. A
section of $L^k_r$ is a continuous section of $L^k_r
\rightarrow M_r$ which lifts to a $G$-invariant $\Ci$ section of
$L^k_{r,U}$ for every $U$. Since 
every $\T^d$-invariant section of $L^k$ restricts to a $G$-invariant
section of $L^k_{r,U}$, the map (\ref{def2}) is well-defined.
Continuing in this way, we can introduce the K{\"a}hler structure of
$M_r$, the Hermitian and holomorphic structures of $L^k_r$, its
connection and verify that we obtain a well-defined map  $V_k$ as in
definition \ref{defGSmap}.

It is also useful to consider $P$ and $P_\C$ as orbifolds and the
projections $p: P \rightarrow M_r$ and $p_{\C} : P_{\C} \rightarrow
M_r$ as orbifold maps. For instance, let $(|U|,U,G,\pi_U)$ be a chart
defined as above. Let $$V := \T^d \times \td \times U, \qquad |V| :=
\T^d_{\C}.U$$ 
and $\pi_V$
be the map $V \rightarrow |V|$ which sends 
$(\te,t,u)$ into $l_{\te+ it}.u$. Let $G$ acts on $V$ by 
$$ G \times V \rightarrow V, \quad g, (\te,t,u) \rightarrow ( \te-g,t, l_{g}.u)$$
Then $(|V|,V,G,\pi_V)$ is an orbifold chart of $P_\C$. Furthermore
$\T^d_\C$ acts on $V = \T_\C^d \times U$ by left multiplication, this action lifts the $\T^d_\C$-action on
$|V|$, and the projection $V \rightarrow U$ locally lifts $p_\C$:
\begin{gather} \label{liftpC}
 \begin{CD}   V @>>> U \\ @VV {\pi_V} V  @VV {\pi_U} V  \\  |V| @>>>
   |U| \end{CD} \quad . \end{gather}
Now, instead of viewing $U$ as a submanifold of $P$, we consider it as the quotient
of $V$ by $\T^d_\C$. To define the
various structure on $U$,  we can lift everything from $|V|$
to $V$ and we perform the reduction from $V$ to $U$. Since the 
$\T^d_\C$-action on $V$ is free, we are reduced to the
manifold case. 
Furthermore since $V \rightarrow
|V|$ is a $G$-principal bundle and $V \rightarrow U$ is
$G$-equivariant, we obtain $G$-invariant structures.
We can apply the same method with the map $p: P \rightarrow
M_r$.
\qed \end{rem}

\begin{rem} (proof of theorem \ref{GuStTh}). Since $P_\C$ is open, $\op{V}_k$ is injective. That $\op{V}_k$ is surjective is more difficult to prove. It consists to show that every invariant holomorphic section of $L^k \rightarrow P_\C$ extends to an invariant holomorphic section over $M$. 
Let us precise that the proof of Guillemin and Sternberg extends to the
orbifold case without modification. The only technical point is to
show that there exists a non-vanishing section in
${\mathcal{H}}_{\la,k}$, when $k$ is sufficiently large (theorem 5.6
of \cite{GuSt1}). This will be proved in section \ref{pilapila}, cf. remark
\ref{comple}.  
\qed \end{rem}

\begin{rem} (harmonic oscillator). The definition of $\op{V}_k :
  {\mathcal{H}}_{1,k} \rightarrow {\mathcal{H}}_{r,k}$ is the
  same. It is easily checked that this map is onto: every
holomorphic section of $L^k_r$ lifts to a holomorphic section $\Psi$
of $L^k$ over $\C^n -\{ 0 \}$ satisfying 
\begin{gather} \label{dsada}
 \Psi(u^{\p_1}z_1,..., u^{\p_n}z_n) = u^k   \Psi(z_1,...,z_n). \end{gather} 
Since it is
bounded on a neighborhood of the origin, it extends on
$\C^n$. Writing its Taylor expansion at the origin, we deduce from
(\ref{dsada}) that $\Psi$ is polynomial and belongs to
${\mathcal{H}}_{1,k}$. \qed
\end{rem}


\section{Reduction of Toeplitz operators} \label{rTo}

\subsection{Toeplitz operators}

Let us denote by
$L^2 (M, L^k)$ the space of $L^2$ sections of $L^k$.
 We define the scalar product of sections of $L^k$ as 
$$ (\Psi,\Psi')_M = \int_{M} (\Psi,\Psi')(x) \ \delta_M (x)$$
where $(\Psi,\Psi')$ is the punctual scalar product and $\delta_M$ is the
Liouville measure $\frac{1}{n!}|\om^{\wedge n}|$. Let $\Pi_k$ be the orthogonal projector of $L^2(M,L^k)$ onto
${\mathcal{H}}_k$.

Given $f \in \Ci(M)$, we denote by $M_f$ the operator of $L^2(M,L^k)$
sending $\Psi$ into $f \Psi$. 
The set of symbols $S(M)$ consists of the sequences $(f(.,k))_k$ of $\Ci(M)$ which admit an asymptotic expansion of the form 
\begin{gather} \label{defs}
 f(.,k) = \sum_{l=0}^{\infty}  k^{-l} f_l +O(k^{-\infty}), \qquad \text{ with } f_0,f_1,.. \in \Ci(M),
\end{gather}
for the $\Ci$ topology.  

A Toeplitz operator is a family $(\op{T}_k)_{k}$ of the form  
\begin{gather} \label{defT} 
 \op{T}_k = \Pi_k M_{f(.,k)} \Pi_k + \op{R}_k \end{gather} 
where $(f(.,k)) \in S(M)$ and $\op{R}_k$ is an operator of
$L^2(M,L^k)$ satisfying $\Pi_k \op{R}_k \Pi_k = \op{R}_k$ and whose
uniform norm is $O(k^{-\infty})$.
The following result is a consequence of the works of Boutet de Monvel
and Guillemin \cite{BoGu} (cf. \cite{oim1}). 
\begin{theo}  \label{ttoep}
The set ${\mathcal{T}}$ of Toeplitz operators is a $*$-algebra. The contravariant symbol map 
$$ \si_{\op{cont}} : {\mathcal{T}} \rightarrow \Ci(M)[[\hb]], \quad \Pi_k M_{f(.,k)} \Pi_k + \op{R}_k \rightarrow \textstyle{\sum} \hb^l f_l $$
is well-defined, onto and its kernel consists of the Toeplitz operators whose uniform norm is $O(k^{-\infty})$. Furthermore, the product $*_{c}$ induced on $\Ci(M)[[\hb]]$ is a star-product. 
\end{theo}

The {\em principal} symbol of a Toeplitz operator is the first coefficient
$f_0$ of its contravariant symbol. The operators $\muquant_{\xi,k}$
are Toeplitz operators with principal symbol  $\langle \mu , \xi \rangle$.

We use the same definitions and notations over $M_r$. Recall that $\la$
is fixed and $k$ runs over the positive multiples of $\ka$. So a
Toeplitz operator of ${\mathcal{T}}_r$ is a family 
$$ (\op{T}_k)_{k=\ka, 2\ka,... }.$$
We denote by $\Pi_{r,k}$ the orthogonal projector onto ${\mathcal{H}}_{r,k}$ and by $*_{cr}$ the product of the contravariant symbols of $\Ci(M_r)[[\hb]]$. 

\begin{rem} (Orbifold) In the case $M_r$ is an orbifold, the
  definition of the Toeplitz operators makes sense. We will prove
  theorem \ref{ttoep} for the Toeplitz operators of $M_r$ in section \ref{ToepOrb}.
 \qed
\end{rem} 

\begin{rem} (harmonic oscillator)
To avoid a discussion about the infinity of $\C^n$, we do not define
the full algebra of Toeplitz operators on $\C^n$ and do not state any result similar to theorem \ref{ttoep}.
We only consider the Toeplitz operators of the form
$$ \Pi_{k} M_{f(.,k)} \Pi_k $$ 
where $f(.,k)$ is a symbol of $S(\C^n)$ (cf. definition in section
\ref{redosc}). \qed
\end{rem}

\subsection{Statement of the main result}

 Recall that $\op{V}_k$ is the isomorphism from ${\mathcal{H}}_{\la,k}$
to  ${\mathcal{H}}_{r,k}$ (cf. definition \ref{defGSmap}). Let $\op{U}_k$
be the operator 
$$ L^2(M, L^k) \rightarrow L^2(M_r, L_r^k), \quad \Psi \rightarrow
\begin{cases} \op{V}_k  ( \op{V}_k^* \op{V}_k)^{-\frac{1}{2}} \Psi,
  \quad \text{
    if } \Psi \in {\mathcal{H}}_{\la,k} \\ 0, \quad  \text{ if $ \Psi$ is
    orthogonal to ${\mathcal{H}}_{\la,k}$} \end{cases} $$ 
Hence 
$$ \op{U}_k^* \op{U}_k =\op{\Pi}_{\la,k}, \quad \op{U}_k \op{U}_k^* =\op{\Pi}_{r,k}, \quad \Pi_{r,k} \op{U}_k \op{\Pi}_{\la,k} = \op{U}_k.
$$
where $\op{\Pi}_{\la,k}$ is the orthogonal projector onto ${\mathcal{H}}_{\la,k}$. 
The main result of the section is the following theorem and the
corresponding theorem
\ref{redToeposc} for the harmonic oscillator.  
\begin{theo} \label{ToepRedToep}
Let $\op{T}_k$ be a Toeplitz operator of $M$ with principal symbol $f$. Then $$ \op{U}_k\op{T}_k\op{U}_k^*: L^2(M_r, L_r^k) \rightarrow  L^2(M_r, L_r^k)$$
 is a Toeplitz operator of $M_r$. Its principal symbol is the function $g \in \Ci(M_r)$ such that 
$$ g \bigl( p(x) \bigr) = \int_{\T^d} f( l_{\te}. x) \; \delta_{\T^d}
(\te ), \quad x \in P, $$
with $\delta_{\T^d}$ the Haar measure of $\T^d$.
\end{theo}

In the following subsection we introduce the $\la$-Toeplitz operators. These are the operators of the form
$$ \Pi_{\la,k} \op{T}_k \Pi_{\la,k} \quad \text{ where } \op{T}_k \in {\mathcal{T}} .$$ 
In the next subsections, we prove some estimates for the sections of
${\mathcal{H}}_{\la,k}$ and introduce an integration map. Then we
prove that the space of $\la$-Toeplitz operators is isomorphic to
the space of Toeplitz operators of $M_r$, a stronger result than theorem
\ref{ToepRedToep}. 

Our proof uses the properties
of the Toeplitz operators of $M_r$ stated in theorem
\ref{ttoep}. So, in the case the $\T^d$-action is not free, the
proof will be complete only in section \ref{ToepOrb} with the proof of
theorem \ref{ttoep} for orbifolds. 

For the following, we introduce the inverse $\op{W}_k :
{\mathcal{H}}_{r,k} \rightarrow  {\mathcal{H}}_{\la,k}$ of $\op{V}_k$. 
We consider that $\op{V}_k$ and $\op{W}_k$ act not only on ${\mathcal{H}}_{\la, k}$ and ${\mathcal{H}}_{r,k}$ respectively, but on the space of $L^2$ sections in such a way that they vanish on the orthogonal of ${\mathcal{H}}_{\la, k}$ and ${\mathcal{H}}_{r,k}$ respectively. So
$$ \Pi_{r,k} \op{V}_k \Pi_{\la,k} = \op{V}_k,\quad   \Pi_{\la,k} \op{W}_k \Pi_{r,k} = \op{W}_k, \quad \op{V}_k \op{W}_k = \Pi_{r,k}, \quad \op{W}_k \op{V}_k = \Pi_{\la,k},$$
and with the convention $0^{-\frac{1}{2}} = 0$, the equality $\op{U}_k = \op{V}_k  ( \op{V}_k^*
\op{V}_k)^{-\frac{1}{2}}$ is valid on $L^{2}(M, L^k)$. 
Furthermore we say that a function or a section is invariant if it
is invariant with respect to the action of $\T^d$. 

\subsection{The $\la$-Toeplitz operators} \label{laTose}

We begin with a useful formula for the orthogonal projector $\laP$ onto $\laH$. Denote by $\op{P}_{\la,k}$ the orthogonal projector of $L^2(M,L^k)$ onto the space of invariant sections of $L^k$ (not necessarily holomorphic). If $\Psi$ is a section of $L^k$, $ \op{P}_{\la,k} \Psi$ is given by the well known formula 
\begin{gather} \label{moyenne}
\op{P}_{\la,k} \Psi =  \int_{\T^d} {\mathcal{L}}_{\te}^* \Psi\;  \de_{\T^d}(\te). \end{gather}
where $\de_{\T^d}$ is the Haar measure. 
Since $\op{P}_{\la,k}$ sends ${\mathcal{H}}_k$ in ${\mathcal{H}}_k$, we have  
\begin{gather} \label{eqes}
 \laP = \op{P}_{\la,k} \Pi_k = \Pi_k \op{P}_{\la,k}. \end{gather} 

Let $\laT$ be the set of $\la$-Toeplitz operators 
$$ \laT := \{ \laP \op{T}_k \laP;\;  \op{T}_k \in {\mathcal{T}}\}.$$
This first result follows from theorem \ref{ttoep} about Toeplitz operators. 

\begin{theo} \label{theolaT}
$ \laT$ is a $*$-algebra. Furthermore if $\op{T}_k$ is a Toeplitz operator with contravariant symbol $\sum \hb^l f_l$, then 
\begin{gather} \label{mueq}
 \laP  \op{T}_k \laP = \laP M_{f_\la(.,k)} \laP + \op{R}_k 
\end{gather}  
where  $\op{R}_k$ is $O(k^{-\infty})$, $ \laP \op{R}_k \laP = \op{R}_k$ and $f_\la(.,k)$ is an invariant symbol of $M$ with an asymptotic expansion $\sum k^{-l} f_{\la,l}$ such that 
$$ f_{\la,l}(x) = \int_{\T^d} f_l(l_{\te}.x) \; \de_{\T^d}(\te), \quad l=0,1,...$$
\end{theo}

So the $\la$-Toeplitz operators can also be defined as the operators of the form 
\begin{gather} \label{LaToEq} 
\laP M_{f(.,k)} \laP + \op{R}_k \end{gather} 
where the multiplicator $f(.,k) \in S(M)$ is invariant, $ \laP \op{R}_k \laP = \op{R}_k$ and $\op{R}_k$ is $O(k^{-\infty})$.

A similar algebra was introduced by Guillemin and Sternberg \cite{GuSt3} in the context of pseudodifferential operators. Their main theorem was that there exists an associated symbolic calculus, where the symbols are defined on the reduced space $M_r$. Let us state the corresponding result in our context. 

\begin{theo} \label{symblaT}
The map $\si_{\op{princ}}: {\mathcal{T}}_\la \rightarrow \Ci(M_r)$, which associates to a $\la$-Toeplitz operator of the form (\ref{LaToEq}) with an invariant multiplicator $f(.,k)$, the function $g_0 \in \Ci(M_r)$ such that 
$$ f(.,k) = p^* g_0 +O (k^{-1})  \quad \text{ over $P$ }$$ 
is well-defined. Furthermore the following sequence is exact
\begin{gather*}  \begin{CD}
 0 @>>> {\mathcal{T}}_\la \cap O(k^{-1}) @>>>  {\mathcal{T}}_\la
 @>\si_{\op{princ}}>> \Ci(M_r) @>>>0 .
\end{CD} \end{gather*} 
Finally if $\op{T}_k^1$ and $\op{T}_k^2$ are $\la$-Toeplitz operators, then
\begin{gather*}  \si_{\op{princ}}( \op{T}_k^1.\op{T}_k^2) = \si_{\op{princ}}( \op{T}_k^1).\si_{\op{princ}}(\op{T}_k^2)
\end{gather*} 
So $ [\op{T}_k^1,\op{T}_k^2]$ is $O(k^{-1})$ and $ k
[\op{T}_k^1,\op{T}_k^2]$ belongs to ${\mathcal{T}}_\la$. Its principal
symbol is 
$$ \si_{\op{princ}}\bigl( k [\op{T}_k^1,\op{T}_k^2]  \bigr) =  i \{ \si_{\op{princ}}( \op{T}_k^1), \si_{\op{princ}}(\op{T}_k^2) \} $$ 
where $\{.,.\}$ is the Poisson bracket of $\Ci(M_r)$. 
\end{theo} 
This theorem doesn't follow from theorem  \ref{ttoep}. Actually it is
a corollary of theorem \ref{tptp}, which says that the algebra of
$\la$-Toeplitz operators and the algebra of Toeplitz operators of
$M_r$ are isomorphic (cf. remark \ref{rr1}). 

\begin{proof}[Proof of theorem \ref{theolaT}]
First let us prove the second point.  Let 
$$\op{T}_k = \Pi_k M_{f(.,k)} \Pi_k + \op{R}_k$$
be a Toeplitz operator. Since $\Pi_k \laP= \laP  \Pi_k = \laP$, we have 
$$  \laP \op{T}_k \laP =  \laP M_{f(.,k)} \laP +  \laP \op{R}_k \laP .$$
Clearly $\op{R}_k$ is $O(k^{-\infty})$ implies that $\laP \op{R}_k \laP$ is $O(k^{-\infty})$.  
Since  $\laP \op{P}_{\la,k} =  \op{P}_{\la,k} \laP =  \laP$, we have 
$$\laP M_{f(.,k)} \laP = \laP \op{P}_{\la,k}  M_{f(.,k)}
\op{P}_{\la,k} \laP.$$ 
Then using (\ref{moyenne}), we obtain 
$$ \op{P}_{\la,k}  M_{f(.,k)} \op{P}_{\la,k} =   M_{f_\la(.,k)} \op{P}_{\la,k}$$ 
where $f_\la(.,k)$ is the invariant symbol
$$ f_\la (x,k) =  \int_{\T^d} f ( l_\te.x,k) \;  \de_{\T^d}(\te). $$
Consequently, 
$$\laP \op{T}_k \laP =  \laP M_{f_\la(.,k)} \laP +  \laP \op{R}_k \laP, $$
which gives the result. To prove that $\laT$ is a $*$-algebra, the only difficulty is to check that the product of two $\la$-Toeplitz operators is a $\la$-Toeplitz operator. Let $f^1(.,k)$ and $f^2(.,k)$ be invariant symbols of $S(M)$. We have to show that
$$ \Pi_{\la,k} M_{f^1(.,k)} \Pi_{\la,k} M_{f^2(.,k)} \Pi_{\la,k} $$ 
is a $\la$-Toeplitz operator. By (\ref{eqes}), 
\begin{align*} 
\Pi_{\la,k} M_{f^1(.,k)} \Pi_{\la,k} M_{f^2(.,k)} \Pi_{\la,k} = & \Pi_{\la,k} M_{f^1(.,k)} \Pi_k \op{P}_{\la,k} M_{f^2(.,k)} \Pi_{\la,k} \\
\intertext{Since $f^2(.,k)$ is invariant, $\op{P}_{\la,k}$ and $M_{f^2(.,k)}$ commute, so}
  = & \Pi_{\la,k} M_{f^1(.,k)} \Pi_k  M_{f^2(.,k)} \op{P}_{\la,k} \Pi_{\la,k} \\
 = & \Pi_{\la,k} \Pi_k M_{f^1(.,k)} \Pi_k  M_{f^2(.,k)} \Pi_k \Pi_{\la,k} \end{align*}
by (\ref{eqes}). Finally $ \Pi_k M_{f^1(.,k)} \Pi_k  M_{f^2(.,k)} \Pi_k$ is a Toeplitz operator since it is the product of two Toeplitz operators. 
\end{proof}

\begin{rem} (harmonic oscillator)
We define the $\la$-Toeplitz operators as the
operators of the form  
$$ \op{\Pi}_{1,k} M_{f(.,k)} \op{\Pi}_{1,k} + \op{R}_k $$
where $f(.,k)$ belongs to $S(\C^n)$, $\op{\Pi}_{1,k}$ is the
orthogonal projector onto ${\mathcal{H}}_{1,k}$,
$\op{R}_k$ satisfies $ \op{\Pi}_{1,k} \op{R}_k
\op{\Pi}_{1,k}=\op{R}_k$ and its uniform
norm is $O(k^{-\infty})$. Then theorems \ref{theolaT}
and \ref{symblaT} remain true. The proof that the multiplicator
$f(.,k)$ can be chosen invariant  is the same. The fact that these
operators form an algebra and the definition and properties of the
principal symbol are consequences of theorem \ref{tptp}. \qed
\end{rem}

\subsection{Norm of the invariant states} \label{secscal}
In this section, we estimate the norm of the eigensections of ${\mathcal{H}}_{\la,k}$. We begin with the estimation over $P_\C$. Using the diffeomorphism (\ref{diffeo}), we identify $P_\C$ with $\td  \times P$. Let $\xi_i$ be a basis of $\td$ and denote by $t_i$ the associated linear coordinates.

\begin{prop} \label{pcontrolnorm}
For every $\T^d_\C$-invariant section $\Psi$ of $L^k \rightarrow P_\C$, we have 
\begin{gather} \label{controlnorm} 
  (\Psi,\Psi) (t,y) = e^{-k \ph(t,y) } ( \Psi, \Psi) (0,y), \quad
  \forall \; (t,y) \in  \td  \times P
\end{gather} 
where $\ph$ is the $\Ci$ function on $\td  \times P$ solution of the equations 
$$ \ph(0,y) = 0, \qquad \partial_{t_i} \ph(t,y) = 2(\la_i - \mu_i(t,y)) , \text{ with } i=1,...,d$$
and $\mu_i := \langle \mu, \xi_i \rangle$,  $\la_i := \langle \la, \xi_i \rangle$ are the components of $\mu$ and $\la$.
\end{prop}

\begin{proof} 
By equation (\ref{repLiecomp}), $i \xi \in i\td$ acts on the sections
of $L^k$ by $$\nabla_{J \xi^\#} - k \langle \mu -\la, \xi\rangle.$$ 
So if $\Psi$ is a $\T^d_\C$-invariant section, 
$$ \nabla_{J \xi^\#}\Psi  = k \langle \mu -\la, \xi \rangle \Psi .$$
which leads to 
\begin{gather} \label{ddcp}
  (J \xi^\#). (\Psi, \Psi) = 2 k \langle \mu -\la, \xi \rangle  (\Psi, \Psi)
\end{gather} 
and shows the proposition. \end{proof}
 
On the complementary set $P^c_\C$ of $P_\C$ in $M$, the situation is simpler. 
\begin{prop} \label{pcompl}
Every $\T^d_\C$-invariant section $\Psi$ of $L^k \rightarrow M$ vanishes over  $P^c_\C$.
\end{prop}
\begin{proof} 
This is also a consequence of equations (\ref{ddcp}) (cf. theorem 5.4 of \cite{GuSt1}). \end{proof}

The previous propositions were shown by Guillemin and Sternberg in \cite{GuSt1}. Furthermore they noticed the following important fact. 

\begin{lemme} \label{derivph}
Let $g$ be the K{\"a}hler metric ($g(X,Y) = \om(X,JY)$). Then 
\begin{gather*} 
 \tfrac{1}{2}\partial_{t_i}\partial_{t_j} \ph(t, y) = g( \xi_i^\#, \xi_j^\# )(t,y). \end{gather*}
\end{lemme} 

\begin{proof} We have $J \xi^\#_i. \langle \mu,\xi_j \rangle = \om ( J \xi^\#_i, \xi_j^\# ) = - g (\xi^\#_i, \xi_j^\# )$. The result follows from proposition \ref{pcontrolnorm}. \end{proof}

Then for every $y \in P$, the function $\ph(.,y)$ is strictly convex. It admits a global minimum at $t =0$ and this minimum is $\ph (0,y) =0$. We obtain the following proposition.

\begin{theo} \label{estimloin}
Let $\eps >0$, $P_\eps$ be the subset of $P_\C$ 
$$ P_\eps := \{ (t,y) \in P_\C; \; |t| < \eps \}$$ 
and $P_\eps^c$ its complementary subset in $M$. There exists some positive constants $C(\eps)$, $C$, $C'$ such that for every $k$ and every $\Psi \in {\mathcal{H}}_{\lambda,k}$,
$$  (\Psi, \Psi) (x) \leqslant C k^{n}  e^{-k C(\eps)} (\Psi,\Psi)_M, \quad \forall x \in P^c_\eps  $$
$$ \text{ and }\quad (\Psi,\Psi) _{P_\eps^c} \leqslant   C' k^{n} e^{-k C(\eps)}
(\Psi,\Psi)_M. $$
where $ (\Psi,\Psi) _{P_\eps^c} := \int_{P_\eps^c}
(\Psi,\Psi) \delta_M$
\end{theo}

\begin{rem}
This result shows that the eigenstates of ${\mathcal{H}}_{\la,k}$ are
concentrated on $P$. Actually, since the $\op{T}_{\xi,k}$ are Toeplitz
operators with principal symbol $\langle \mu, \xi \rangle$, we could
directly deduce from general properties of these operators \cite{oim2} a weaker version where $C k^{n}  e^{-k C(\eps)}$ is replaced by $C_N(\eps) k^{-N}$ with $N$ arbitrary large.  
\qed \end{rem}

\begin{proof} 
Let $C(\eps)$ be the minimum value of $\ph$ over the compact set
$$\{(t,y) \in P_\C; \; |t| = \eps \}.$$ Since $\ph(.,y)$ is strictly
convex with a global vanishing minimum at $t=0$,  $C(\eps)$ is positive and 
\begin{gather*} \label{fgh} 
\ph (t,y) \geqslant C(\eps), \quad \forall \; (y,t) \in P_\eps^c. \end{gather*}
On the other hand, using coherent states as in section 5 of \cite{oim1}, we prove that there exists a constant $C$ such that for every $k$ and $\Psi \in {\mathcal{H}}_k$,
$$ (\Psi, \Psi) (x) \leqslant C k^{n} (\Psi,\Psi)_M,\quad  \forall x \in P.$$
If $\Psi$ belongs to ${\mathcal{H}}_{\la,k}$, then it is $\T^d_\C$-invariant and so it satisfies equation (\ref{controlnorm}). Furthermore it vanishes over the complementary set of $P_\C$. This implies the first part of the result. 
By integrating, we get the second part with $C' = C \op{Vol} (P_\eps)$. 
\end{proof} 

\begin{rem} (harmonic oscillator)
We have to adapt theorem \ref{estimloin} since $\C^n$ is not compact. The imaginary part of the complex action is given by 
$$ \R \times \C^n \rightarrow \C^n, \quad (t,z) \rightarrow  l_{it}.z = ( e^{-t \p_1}z_1,..., e^{-t \p_n}z_n).$$ 
A section $\Psi$ of ${\mathcal{H}}_{1,k}$ satisfies  
\begin{gather} \label{lkj}
 \Psi(l_{it}. z) = e^{-kt} \Psi(z) \end{gather}  
From this we obtain the following lemma. 

\begin{lemme} \label{estimCn}
Let $g(.,k)$ be a sequence of functions on $\C^n$. Assume that there exists $C$ and $N$ such that 
$$ | g(z,k)| \leqslant C(1+ |z| )^{N}, \quad \forall z \in \C^n, \;\forall k. $$ 
Let $\eps >0$, $P_\eps$ be the subset $\{ l_{it}.z; \; |t| < \eps \text{ and } z \in P\}$ of $P_\C$ and $P_\eps^c$ its complementary subset in $\C^n$. 
There exists some positive constants $C(\eps)$, $C'$ such that for every $k$ and every $\Psi \in {\mathcal{H}}_{1,k}$,
$$  \quad ( g(.,k) \Psi,\Psi) _{P_\eps^c} \leqslant   C' k^{n} e^{-k C(\eps)} (\Psi,\Psi)_{\C^n}. $$
\end{lemme}
\begin{proof} 
Let $\Psi$ belongs to ${\mathcal{H}}_{1,k}$. First, as in the compact case, there exists $C_1$ such that 
$$ |\Psi(z)|^2 e^{-k |z|^2} \leqslant C_1 k^{n} (\Psi, \Psi)_{\C^n}, \quad \forall\; z \in P$$ 
Let $w = l_{it}.z$ with $z \in P$. It follows from (\ref{lkj}) that 
\begin{xalignat*}{2}  
|\Psi( w) |^2 e^{-k | w | ^2} = & |\Psi( z) |^2 e^{-k | z| ^2} e^{ -k
 ( |w|^2  + 2t - |z| ^2  )} \\
\leqslant &   C_1 k^{n}  e^{ -k
 ( |w|^2  + 2t - |z| ^2  )}(\Psi, \Psi)_{\C^n}
\end{xalignat*} 
Using that $|z|^2 \leqslant 1$, we obtain
\begin{gather}  \label{lala}  
|\Psi( w) |^2 e^{-k | w | ^2} \leqslant   C_1 k^{n}  e^{ -k }
(\Psi, \Psi)_{\C^n}\qquad \text{ if $t \geqslant 1$.} \end{gather} 
On the other hand, assume that $t <0$. There exists $C_2 >0$ such that
$|z|^2 \geqslant C_2$ for
every $z \in P$.  So $| w|^2 = \textstyle{\sum} e^{-2 \p_i t} |z_i|^2
\geqslant C_2 e^{-2t} $. Consequently 
 \begin{xalignat*}{2} e^{ -k
 ( |w|^2  + 2t - |z| ^2  )} \leqslant & C_3  e^{- k( |w|^2 -2 \ln |w|
 -1 )}\\ \leqslant & C_3  e^{- \frac{k}{2} |w|^2} e^{-k}, \qquad \text{ if $|w|$ is sufficiently large.} \end{xalignat*} 
Hence there exists $t_- <0$ such that 
$$ |\Psi( w) |^2 e^{-k | w | ^2} \leqslant  C_4  k^{n}  e^{ -k }
(\Psi, \Psi)_{\C^n} e^{- \frac{k}{2} |w|^2} \qquad \text{ if $t \leqslant t_-$.}
$$ 
This last equation and \eqref{lala}  lead to the result if
$\epsilon \geqslant \max (1, - t_-)$. For the smaller values of $\epsilon$, we complete the proof as in theorem \ref{estimloin} since we are reduced to a compact subset of $P_{\C}$. 
\end{proof}
\end{rem}

\subsection{The integration map} \label{secIM}

Recall that $\delta_M$ and $\delta_{M_r}$ are the Liouville measures of $M$ and $M_r$ respectively and $p_{\C}$ denote the projection $P_{\C} \rightarrow M_r$.

 Let $I_k : \Ci_{o}( P_\C) \rightarrow \Ci(M_r)$ be the map given by 
$$  I_k(f) (x) = \Bigl( \frac{k}{2\pi} \Bigr)^{\frac{d}{2}} \Biggr( \int_{ p_{\C}^{-1}(x)} e^{-k \ph} f \delta_M \Biggl). \delta_{M_r}^{-1}(x) $$
where $\ph$ is defined in proposition
\ref{pcontrolnorm}. Equivalently, 
$I_k(f) \; \delta_{M_r}$ is the
push-forward of $(k/2\pi)^{\frac{d}{2}}e^{-k \ph} f \; \delta_M$ by $p_\C$. 

\begin{rem} \label{Orbor} (Orbifold). If $M_r$ is an orbifold, the
  push-forward by $p_\C$ of a density $\nu$ of $P_\C$  with compact
  support is defined as in the manifold case in such a way that:
$$ \forall g \in \Ci(M_r), \int_{P_\C} \nu\; p_\C^*g =  \int_{M_r}  g\;
p_{\C *} \nu. $$
Applying this in orbifold charts of $P_\C$ and $M_r$, we recover the
usual definition. With the same
notations as in remark \ref{orbidata}, assume that $g$
has a compact support in $|U|$. Denote by $g_U$, $(p_{\C *} \nu)_U$
the local lifts in $U$ of $g$ and $p_{\C *} \nu$. Then by definition of an integral in an orbifold,
$$ \int_{M_r}  g\; p_{\C *} \nu = \frac{1}{\# G} \int_U
g_U. (p_{\C *} \nu)_U $$
So it follows from (\ref{liftpC}) that $(p_{\C *} \nu)_U$ is the
push-forward of $\pi_V^*\nu$ by the projection $V \rightarrow U$. One
can check that the $(p_{\C *} \nu)_U$ agree on overlaps and define a
global section $p_{\C *}  \nu$. \qed
\end{rem}

\begin{rem} (harmonic oscillator). Since we apply $I_k $ only to functions
  with compact support, all the results in this section extend directly to
  this case. \qed
\end{rem}

We introduced the map $I_k$ because it satisfies the following
property. 
\begin{prop} \label{pres}
For every $f \in \Ci_{o}( P_\C)$, we have 
\begin{gather} \label{propI} 
 ( f \Psi, \Psi')_M = ( I_k(f)  \op{V}_k \Psi, \op{V}_k \Psi' )_{M_r} ,\qquad  \forall \Psi, \Psi' \in {\mathcal{H}}_{\la,k} \end{gather} 
or equivalently 
$$ \Pi_{\la,k} M_{f} \Pi_{\la,k} = \op{V}_k^* M_{I_k(f)} \op{V}_k. $$
\end{prop}
\begin{proof}  Let us prove (\ref{propI}). Since the support of $f$ is a subset of $P_\C$, 
\begin{align*} 
( f \Psi, \Psi')_M =&  \int_{P_\C} f \; (  \Psi, \Psi') \;  \delta_{M}
\end{align*} 
By definition of
$\op{V}_k$, we have 
$$   p^* (\op{V}_k \Psi, \op{V}_k \Psi') = \Bigl( \frac{2\pi}{k}
\Bigr)^{\frac{d}{2}} j^* (\Psi,\Psi')  . $$
So equation (\ref{controlnorm}) can be rewritten as 
$$ (  \Psi, \Psi') =  \Bigl( \frac{k}{2\pi} \Bigr)^{\frac{d}{2}} e^{-k \ph} p_{\C}^* (\op{V}_k \Psi, \op{V}_k \Psi').$$ 
Hence, 
\begin{align*} 
( f \Psi, \Psi')_M = & \Bigl( \frac{k}{2\pi} \Bigr)^{\frac{d}{2}}
\int_{P_\C} f \; e^{-k \ph} p_{\C}^* (\op{V}_k \Psi, \op{V}_k \Psi') \\
= &  \int_{M_r} I_k(f) \; (\op{V}_k \Psi, \op{V}_k
\Psi') \;  \delta_{M_r}\\
= & ( I_k(f)  \op{V}_k \Psi, \op{V}_k \Psi' )_{M_r}.
\end{align*} 
which proves the result.   
\end{proof}
For the following, we need to control the asymptotic behavior of $ I_k( f)$ as $k$ tends to infinity when $f$ depends also on $k$. First observe that there is no restriction to consider only invariant functions of $\Ci_{o}( P_\C)$. Indeed, if $f\in \Ci_{o}( P_\C)$ and 
$$ f_\la (x ) = \int_{\T^d} f (l_\te.x) \; \de_{\T^d}(\te) $$
then $I_k(f) = I_k( f_\la)$. 

Denote by $S_o(P_\C)$ the set of sequences $(f(.,k) )_k$ of $\Ci_{\T^d}(P_\C)$ such that there exists a compact set $K \subset P_\C$ which contains the support of $f(.,k)$ for every $k$ and $f(.,k) $ admits an asymptotic expansion for the $\Ci$ topology of the form 
$$ f(.,k) = \textstyle{\sum}_l k^{-l} f_l + O(k^{-\infty}).$$  

Let us introduce a basis $\xi_i$ of the integral lattice of
$\td$ and denote by $\xi^{\#}_i$ the associated vector fields of $M$. Recall that $g$ is the K{\"a}hler metric. The following result
involves the determinant of $ g(\xi_i^\#, \xi_j^\#)$, which clearly
doesn't depend on the choice of the basis $\xi_i$.

\begin{prop} \label{devasI} 
Let $f \in S_o(P_\C)$. Then the sequence $I_k(f(.,k))$ is a symbol of $M_r$. Furthermore 
 $$p^* g_0  =  j^* \bigl( \op{det} \bigl[ g(\xi_i^\#, \xi_j^\#) \bigr]^{\frac{1}{2}} . f_0 \bigr) $$
where $g_0$ and $f_0$ are such that  $ I_k(f(.,k)) = g_0 + O(k^{-1})$ and  $ f(.,k) = f_0 + O(k^{-1})$.
\end{prop}

This proposition admits the following converse. 
\begin{cor} \label{convI}
For every symbol $g(,.k)$ of $S(M_r)$, there exists $f(.,k) \in S_o(P_\C)$ such that $g(.,k) = I_k(f(.,k))$. \end{cor}
\begin{proof}[Proof of corollary \ref{convI}]
Let $r$ be a non negative invariant function of $\Ci_o(P_\C)$ such that $r = 1$ on a neighborhood of $P$. Set 
$$ f (.,k) :=  r. p_{\C}^*\bigl( g(.,k) I^{-1}_k(r) \bigr) .$$ 
From the previous proposition, $I_k(r)$ is a symbol of $M_r$ and the
first coefficient of its asymptotic expansion doesn't vanish. So
$I^{-1}_k(r)$ is also a symbol of $M_r$. Consequently $ f (.,k) $
belongs to $S_o(P_\C)$. Furthermore, we have $g(.,k) = I_k(f(.,k))$. 
\end{proof}

\begin{proof}[Proof of proposition \ref{devasI}]  
We integrate first over the fibers of $\td \times P \rightarrow
P$. Let us compute the Liouville measure $\de_M$. We denote by
$t_1,...,t_d$ the linear coordinates of $\td$ associated to
$\xi_1,...,\xi_d$.   
\begin{lemme} There exists an invariant measure $\delta_P$ on $P$
  and a function $\delta \in \Ci(\td \times P)$ such that
\begin{gather*} 
 \delta_M = \delta. \delta_{P}.|dt_1...dt_d|  \text{ over } \td \times
 P, \\
 p_* \delta_P = \delta_{M_r} \quad \text{ and } \quad \delta (0,.) = j^* \det [ g(\xi_i^\#, \xi_j^\#) ].
\end{gather*}  
\end{lemme} 
\begin{proof}
Let us write over $\{ 0 \} \times P \subset \td \times P$
$$ \om = \be + \sum_{ 1 \leqslant i\leqslant d} \be_i \wedge dt_i + \sum_{1 \leqslant i<j\leqslant d} a_{ij} dt_i \wedge dt_j $$ 
where $\be \in \Om^2(P)$, $\be_i \in \Om^1(P)$ and $a_{ij} \in
\Ci(P)$. Since this decomposition is unique, these forms are all
invariant. 
Let  us set 
$$ \delta_P =  
\frac{ | \be^{ \wedge n_r}|}{n_r !} . \frac{ |\be_1 \wedge ...\wedge \be_d|} {\op{det} [ g(\xi_i^\#, \xi_j^\#) ]} $$
Since $j^* \om = \be$, $\be = p^*\om_r$. We also have $g(\xi_i^\#,
\xi_j^\#)= \om ( \xi_i^\#, J \xi_j^\#) = \langle \be_j, \xi_i^\#
\rangle$. Since the $\xi_i$ are a basis of the integer lattice, we
obtain $p_* \delta_P = \delta_{M_r}$. In the case $M_r$ is an
orbifold, this can be proved using local charts of $M_r$ and $P$ as in remark \ref{Orbor}. 

Since $\be (\xi_i^{\#},.)= 0$, $(dt_i \wedge dt_j) (\xi_k^\#, .) =0$ and
  $ (\be_i \wedge dt_i) ( \xi_k^{\#},\xi_l^{\#}) =0$, we  have over $\{ 0 \} \times P$ 
\begin{align*} 
 \om^{\wedge n}
  = & \frac{n!}{n_r !} \be^{ \wedge n_r}\wedge (\be_1 \wedge dt_1) \wedge...\wedge ( \be_d \wedge dt_d )
\end{align*} 
Hence
$$ \delta_M =  \det [ g(\xi_i^\#, \xi_j^\#) ] \; \delta_{P}.|dt_1...dt_d| $$
over $\{ 0 \} \times P$, which proves the result.\end{proof}

Let $J_k(f)$ be the function of $\Ci(P)$ 
$$ J_k(f)(y) = \Bigl( \frac{k}{2\pi} \Bigr)^{\frac{d}{2}} \int_{\td} e^{-k \ph(t,y) } f(t,y) \delta (t,y) \;  |dt_1...dt_d|$$ 
It is invariant and $p^* I_k(f) = J_k(f)$. So we just have to estimate $ J_k(f)(y)$ which can be done with the stationary phase lemma. Recall that we computed in lemma \ref{derivph} the second derivatives of $\ph$. The result follows.\end{proof}

\begin{rem}
The proof actually gives more about the map 
\begin{gather} \label{definF}
 F : \Ci_{\T^d} (M) [[\hb]] \rightarrow \Ci(M_r) [[\hb]], \quad
 \textstyle{\sum} \hb^l f_l \rightarrow  \textstyle{\sum} \hb^l g_l
\end{gather} 
such that $I_k(f(.,k)  ) =\sum k^{-l} g_l +O(k^{-\infty})$ if $f(.,k) = \sum k^{-l} f_l +O(k^{-\infty})$.

Since $F$ enters in the computation of the contravariant symbol of the
reduced operator, let us give its properties. First it is $\C[[\hb]]$-linear. So 
$$ F = \textstyle{\sum} \hb^l F_l \quad \text{with } F_l: \Ci_{\T^d} (M)\rightarrow \Ci(M_r). $$
The operators $F_l$ are of the following form 
$$ p^* F_l(g) =  \op{det}^{\frac{1}{2}} [g(\xi^\#_i ,\xi^\#_j )]
\sum_{|\al| \leqslant 2l}  a_{\al,l} \;  j^* \Bigl( (J\xi_1^\#)^{\al(1)} ...(J\xi_d^\#)^{\al(d)}. f   \Bigr) $$
where the functions $a_{l,\al}$ are polynomials in the derivatives of
$\mu_i =\langle \mu,\xi_i \rangle$ and $\Delta \mu_i$ with respect to
the gradient vector fields of the $\mu_j$. 

Indeed by proposition \ref{pcontrolnorm} the derivatives of $\ph$ can be computed in terms of the derivatives   of $\mu_i$. Furthermore it is easily proved that 
$$ (J\xi_i^\# ).\ln  \de =  \Delta \mu_i $$
with $\Delta$ the Laplace-Beltrami operator of $M$, which gives the
derivatives of $\de$ in terms of the derivatives of $\Delta
\mu_i$. Then the computation of the functions $a_{\al,l}$ follows from
the stationary phase lemma.
\qed \end{rem}

\subsection{From the $\la$-Toeplitz operators to the reduced Toeplitz
  operators} \label{secTT}

We begin with a rough estimate of the maps $\op{V}_k$ and $\op{W}_k$.
\begin{prop} \label{rough} 
There exists a constant $C >0$, such that for every $k$ the uniform norms of $\op{V}_k$, $\op{V}^*_k$, $\op{W}_k$ and $\op{W}_k^*$ are bounded by $C$. 
\end{prop}

\begin{proof} 
By corollary \ref{convI}, there exists $f(,.k) \in S_o(P_\C)$
such that $I_k (f(,.k)) =1$. By proposition \ref{pres},  
$$ \op{V}_k^* \op{V}_k = \Pi_{\la,k} M_{f(.,k)} \Pi_{\la,k}. $$
Furthermore it follows from proposition \ref{devasI} that $f(.,k) = f_0 + O(k^{-1})$ with $f_0$ positive on $P$.  

Since $f(.,k)$ is a symbol, there exists $C_1$ such that $f(.,k) \leqslant C_1$ over $M$ for every $k$. So 
$$ ( \op{V}_k \Psi,  \op{V}_k \Psi)_{M_r} = ( f(.,k) \Psi,   \Psi)_{M} \leqslant  C_1(\Psi, \Psi)_M $$
which proves that the uniform norms of $\op{V}_k$ and  $\op{V}^*_k$ are
smaller than $C_1^{\frac{1}{2}}$. 

Since $ f_0$ is positive on $P$, there exists a neighborhood $P_\epsilon$ of $P$ defined as in theorem \ref{estimloin} and a constant $C_2>0$, such that $$ f(.,k) \geqslant  C_2, \quad \text{ over } P_\epsilon $$
when $k$ is sufficiently large. So 
$$ ( \op{V}_k \Psi,  \op{V}_k \Psi)_{M_r} = ( f(.,k) \Psi,   \Psi)_{M} \geqslant   ( f(.,k) \Psi,   \Psi)_{P_\eps} \geqslant C_2(\Psi, \Psi)_{P_\epsilon} $$
Furthermore theorem \ref{estimloin} implies that 
$$ (\Psi, \Psi)_{P_\epsilon} = (\Psi, \Psi)_{M} - (\Psi, \Psi)_{P_\epsilon^c} \geqslant \tfrac{1}{2} (\Psi, \Psi)_{M}$$ 
when $k$ is sufficiently large. Consequently  the uniform norms of
$\op{W}_k$ and  $\op{W}^*_k$ are smaller than $(\frac{1}{2}C_2)^{-\frac{1}{2}}$ when $k$ is sufficiently large. 
\end{proof}
\begin{rem} (harmonic oscillator) The result is still valid. There are some modifications in the
  proof. Instead of theorem \ref{estimloin}, we have to use lemma
  \ref{estimCn}. The same holds for theorems \ref{keytheo} and \ref{tptp}. \qed \end{rem}

Let us now give the relations between the $\la$-Toeplitz operators and
the Toeplitz operators of $M_r$. 

\begin{theo} \label{keytheo}
If $ \op{T}_k$ is a $\la$-Toeplitz operator, then $\op{W}_k^* \op{T}_k \op{W}_k$ is a Toeplitz operator of $M_r$. Furthermore, if 
$$\op{T}_k = \laP M_{f(.,k)} \laP + O(k^{-\infty}),$$ with $f(.,k) = \sum k^{-l} f_l + O(k^{\infty})$ an invariant symbol, then 
$$ \si_{\op{cont} } ( \op{W}_k^* \op{T}_k \op{W}_k ) = F(
\textstyle{\sum} \hb^l f_l )$$
where the map $F$ is defined in (\ref{definF}).   
Conversely, if $ \op{T}_k$ is a Toeplitz operator of $M_r$, then $\op{V}_k^* \op{T}_k \op{V}_k$ is a $\la$-Toeplitz operator. So the map 
$$ \laT \rightarrow \rT, \quad \op{T}_k \rightarrow \op{W}_k^* \op{T}_k \op{W}_k$$ 
is a bijection. 
\end{theo}

\begin{proof} 
Let 
$$\op{T}_k =  \laP M_{f(.,k)} \laP + \op{R}_k$$
be a $\la$-Toeplitz operator, where $f(.,k) = \sum k^{-l} f_l + O(k^{\infty})$ is an invariant symbol and $\op{R}_k$ is $O(k^{-\infty})$. Then 
$$ \op{W}_k^* \op{T}_k \op{W}_k = \op{W}_k^* M_{f(.,k)} \op{W}_k +  \op{W}_k^*\op{R}_k \op{W}_k $$
By proposition \ref{rough}, $\op{W}_k^*\op{R}_k \op{W}_k$ is $O(k^{-\infty})$. 

Let $P_\epsilon$ be a neighborhood of $P$ defined as in theorem \ref{estimloin}. Let $r$ be an invariant function of $\Ci_o (P_\C)$ such that $r =1$ over $P_\epsilon$. Write 
$$\op{W}_k^* M_{f(.,k)} \op{W}_k = \op{W}_k^* M_{r f(.,k)} \op{W}_k +
\op{W}_k^* M_{(1-r) f(.,k)} \op{W}_k $$

The second term on the right side is $O(k^{-\infty})$. Indeed by proposition \ref{rough}, it suffices to prove that $\laP M_{(1-r) f(.,k)}\laP$ is $O(k^{-\infty})$. We have 
\begin{align*} 
( (1-r) f(.,k) \Psi, (1-r) f(.,k) \Psi )_M & = ( (1-r) f(.,k) \Psi, (1-r) f(.,k) \Psi )_{P^{c}_\epsilon} 
\intertext{since $r =1$ over $P_\epsilon$}
& \leqslant C ( \Psi,  \Psi )_{P^{c}_\epsilon} \end{align*} 
where $C$ doesn't depend of $k$. Theorem \ref{estimloin} leads to the conclusion. 

Now $r f(.,k)$ is a symbol of $S_o(P_\C)$. So by proposition
\ref{devasI}, $g(.,k) := I_k ( r f(.,k))$ is a symbol of $S(M_r)$. By
proposition \ref{pres},  
$$ \op{W}^*_k M_{r f(.,k)} \op{W}_k = \rP M_{g(.,k)} \rP.$$
which proves that $\op{W}^*_k \op{T}_k \op{W}_k$ is a  Toeplitz operator of $M_r$ with contravariant symbol $F( \textstyle{\sum} \hb^l f_l ). $

Conversely, let 
$$ \op{T}_k =  \Pi_{r,k} M_{g(.,k)} \Pi_{r,k} + \op{R}_k $$ 
be a Toeplitz operator of $M_r$, where $ g(.,k) \in S(M_r)$ and $ \op{R}_k$ is $O(k^{-\infty})$. Write 
$$ \op{V}^*_k \op{T}_k \op{V}_k =  \op{V}^*_k M_{g(.,k)} \op{V}_k +  \op{V}_k^* \op{R}_k \op{V}_k . $$
Then, $ \op{V}_k^* \op{R}_k \op{V}_k$ is $O(k^{-\infty})$ by proposition \ref{rough}. By corollary \ref{convI}, there exists $f(.,k) \in S_o(P_\C)$ such that 
$$ \op{V}_k^* M_{g(.,k)} \op{V}_k = \laP M_{f(.,k)} \laP .$$ 
Consequently $\op{V}^*_k \op{T}_k \op{V}_k$ is a $\la$-Toeplitz operator. 
\end{proof}

Recall that $\op{U}_k = \op{V}_k  ( \op{V}_k^*
\op{V}_k)^{-\frac{1}{2}} :  L^2(M, L^k) \rightarrow L^2(M_r,
L_r^k)$. Let us state our main result. 

\begin{theo} \label{tptp}
The map 
 $$ \laT \rightarrow \rT, \quad \op{T}_k \rightarrow \op{U}_k \op{T}_k \op{U}_k^*$$ 
is an isomorphism of $*$-algebra. Furthermore, 
if 
$$\op{T}_k = \laP M_{f(.,k)} \laP + O(k^{-\infty}),$$ with $f(.,k) = \sum k^{-l} f_l + O(k^{\infty})$ an invariant symbol, then 
$$ \si_{\op{cont}} ( \op{U}_k \op{T}_k \op{U}_k^* ) = e^{-\frac{1}{2}} *_{\op{cr}} F( \textstyle{\sum} \hb^l f_l ) *_{\op{cr}}  e^{-\frac{1}{2}} $$  
where $e =  F( 1)$ and $e^{-\frac{1}{2}}$ is the formal series of
$\Ci(M_r)[[\hb]]$ whose first coefficient is positive and such that 
$ e^{-\frac{1}{2}} *_{\op{cr}} e^{-\frac{1}{2}}  *_{\op{cr}} e = 1$.
\end{theo}

Again, in the proof we use some basic properties of the Toeplitz
operators of $M_r$, which are known in the manifold case and will be
extended to the orbifold case in section \ref{ToepOrb}.

\begin{proof}  
It is easily checked that
$$\op{U}_k= ( \op{W}_k^* \op{W}_k ) ^{-\frac{1}{2}} \op{W}_k^*.$$ 
Let $\op{T}_k$ be a $\la$-Toeplitz operator. We have 
\begin{gather} \label{eq1}
 \op{U}_k \op{T}_k \op{U}_k^* =  ( \op{W}_k^* \op{W}_k ) ^{-\frac{1}{2}} \op{W}_k^* \op{T}_k \op{W}_k ( \op{W}_k^* \op{W}_k ) ^{-\frac{1}{2}} \end{gather}
By theorem \ref{keytheo}, $\op{W}_k^* \op{W}_k$ is a Toeplitz operator
of $M_r$ with a positive principal symbol. It follows from the
functional calculus for Toeplitz operator (cf. \cite{oim1}) that $( \op{W}_k^* \op{W}_k ) ^{-\frac{1}{2}}$ is a Toeplitz operator also. Now $\op{W}_k^* \op{T}_k \op{W}_k$ is a Toeplitz operator by theorem \ref{keytheo}. Since the Toeplitz operators of $M_r$ form an algebra, $\op{U}_k \op{T}_k \op{U}_k^*$ is a Toeplitz operator. The computation of its covariant symbol is also a consequence of (\ref{eq1}). Indeed by theorem \ref{keytheo}, the symbol of $\op{W}_k^* \op{T}_k \op{W}_k$ and $\op{W}_k^* \op{W}_k$ are $F( \sum \hb^l f_l)$ and $e$ respectively. 

Conversely, if $\op{S}_k$ is a Toeplitz operator of $M_r$, then 
\begin{align*} 
 \op{U}_k^* \op{S}_k \op{U}_k = &\op{W}_k  ( \op{W}_k^* \op{W}_k )
 ^{-\frac{1}{2}} \op{S}_k ( \op{W}_k^* \op{W}_k ) ^{-\frac{1}{2}}
 \op{W}_k^* \\ = &
(\op{V}_k^* \op{W}_k^*) \op{W}_k  ( \op{W}_k^* \op{W}_k ) ^{-\frac{1}{2}} \op{S}_k ( \op{W}_k^* \op{W}_k ) ^{-\frac{1}{2}} \op{W}_k^* (\op{W}_k \op{V}_k) \\
 = & \op{V}_k^*  ( \op{W}_k^* \op{W}_k ) ^{\frac{1}{2}}  \op{S}_k ( \op{W}_k^* \op{W}_k ) ^{\frac{1}{2}} \op{V}_k.
\end{align*} 
And in a similar way, we deduce from theorem \ref{keytheo} that $ \op{U}_k^* \op{S}_k \op{U}_k$ is a $\la$-Toeplitz operator. 
\end{proof}

\begin{rem} (symbolic calculus). 
Recall that we computed the operator $F$ at the end of section
\ref{secIM}. Furthermore the star-product $ *_{\op{cr}}$ can be
computed in terms of the K{\"a}hler metric of $M_r$
(cf. \cite{oim1}). This leads to the computation of the contravariant
symbol of the reduced operator $\op{U}_k \op{T}_k \op{U}_k^*$ in terms
of the multiplicator $\sum k^{-l} f_l$ defining the $\la$-Toeplitz
operator $\op{T}_k$. In particular, the principal symbol $g_0$ of
$\op{U}_k \op{T}_k \op{U}_k^*$ is such that $ p^* g_0 = i^* f_0$.
\qed \end{rem}

\begin{rem} \label{rr1} (proof of theorem \ref{symblaT}). 
Because of the previous  remark, the $\la$-Toeplitz operator
$\op{T}_k$ and the Toeplitz operator $\op{U}_k \op{T}_k \op{U}_k^*$
have the same principal symbol. 
 Consequently all
the assertions of theorem \ref{symblaT} follow from theorem \ref{tptp} and the
calculus of the contravariant symbol for the Toeplitz
operators of $M_r$. 
\qed \end{rem}

\begin{rem} (proof of theorem \ref{ToepRedToep}).
Theorem \ref{ToepRedToep} stated in the introduction is a consequence of theorems \ref{theolaT} and \ref{tptp}. To compute the contravariant symbol of the reduced operator, we have first to average the contravariant symbol of the Toeplitz operator of $M$ and then apply the formula of theorem \ref{tptp}.
\qed \end{rem}

\begin{rem} ($\op{V}_k$ is not unitary). 
The fact that the Guillemin-Sternberg isomorphism is not unitary, even after the rescaling with the factor $( \frac{k}{2\pi})^{\frac{d}{4}}$, can be deduced from the spectral properties of the Toeplitz operators. Indeed by theorem \ref{keytheo}, $\op{W}_k^*\op{W}_k$ is a Toeplitz operator with principal symbol $g_0$ such that 
$$ p^* g_0 = j^* \op{det}[g(\xi^\#_i,\xi^\#_j)]^{\frac{1}{2}}.$$ 
Denote by $m$ and $M$ the minimum and maximum of $g_0$. Then the smallest eigenvalue $E_s$ of $\op{W}^*\op{W}$ and the biggest $E_S$ are estimated by 
$$ E_s = m + O(k^{-1}), \quad E_S = M + O(k^{-1}). $$  
So when the function $\op{det}[g(\xi^\#_i,\xi^\#_j)]^{\frac{1}{2}}$ is
not constant over $P$, $\op{W}_k$ is not unitary when $k$ is
sufficiently large. This happens for instance in the case of the
harmonic oscillator when the reduced space is not a manifold.
\qed \end{rem}

\begin{rem}
From a semi-classical point of view, the operator $\op{U}_k$ is not unique. Indeed we can replace it with any operator of the form 
$$    \op{T}_k \op{U}_k \op{S}_k$$ 
with $\op{S}_k$ a unitary $\la$-Toeplitz operator and $\op{T}_k$  a unitary Toeplitz operator of $M_r$. We can state a theorem similar to theorem \ref{tptp} with this operator. The only changes are in the symbolic calculus.
\qed \end{rem}


\section{Fourier integral operators} \label{vasavoir}

In this section, we prove that the $\la$-Toeplitz operators, the Guillemin-Sternberg isomorphism and its unitarization are Fourier integral operators in the sense of \cite{oim2}. Using this we can interpret the relations between these operators and the Toeplitz operators as compositions of Fourier integral operators corresponding to compositions of canonical relations.

We assume that the reduced space $M_r$ is a manifold. A part of the
material will be adapted to orbifolds in section \ref{prep}. 

\subsection{Definitions} \label{defFIO}
We first recall some definitions of \cite{oim2}.
Let $M_1$ and $M_2$ be compact K{\"a}hler manifolds endowed with
prequantization bundles $L^{\ka}_1 \rightarrow M_1$ and $L_2^{\kappa}
\rightarrow M_2$. Here $\ka$ is some fixed positive integer and in the
following $k$ is always a positive multiple of $k$. Denote by
${\mathcal{H}}_k^1$ (resp. ${\mathcal{H}}_k^2$) the space of holomorphic
sections of $L^{k}_1$ (resp. $L^k_2$) and by $\Pi^1_k$
(resp. $\Pi^2_k$) the orthogonal projector onto ${\mathcal{H}}_k^1$
(resp. ${\mathcal{H}}_k^2$).

Consider a sequence $(\op{T}_k)_{k\in \ka \N^*}$ such that for every $k$, $\op{T}_k$ is an operator ${\mathcal{H}}_k^2 \rightarrow {\mathcal{H}}_k^1$. As previously we extend $\op{T}_k$ to the Hilbert space of sections with finite norm in such a way that it vanishes on the orthogonal of ${\mathcal{H}}_k^2$. The Schwartz kernel $T_k$ is the section of $L_1^k \boxtimes L_2^{-k} \rightarrow M_1 \times M_2$ such that   
$$ \op{T}_k.\Psi  (x_1) = \int_{M_2} T_k(x_1,x_2). \Psi(x_2) \de_{M_2} (x_2) $$
where $ \de_{M_2}$ is the Liouville measure of $M_2$. All the
operators we consider in this section are of this form.
 
\subsubsection{Smoothing operators} \label{smooth}
A sequence $(f(.,k))$ of functions
on a manifold $X$ is $O_{\infty}(k^{-\infty})$ if for every compact
set $K$, every $N \geqslant 0$, every vector fields $Y_1$,...,$Y_N$ on
$X$ and every $l$, there exists $C$ such that $$   | Y_1.Y_2...Y_N
f(.,k)|  \leqslant C k^{-l} \quad \text{ on $K$}.$$ 

Let $L_X \rightarrow X$ be a Hermitian line bundle. Let $(\Psi_k)$ be a sequence such that for every $k$, $\Psi_k$ is a section of $L^k_X$. Then $(\Psi_k)$ is $O_{\infty}(k^{-\infty})$ if for every local unitary section $t :V \rightarrow L_X$, the sequence $(f(.,k))$ such that $\Psi_k = f(.,k) t^k$ is $O_{\infty}(k^{-\infty})$. 

We say that an operator $(\op{T}_k)$ is smoothing if the sequence
$(T_k)$ of Schwartz kernels is $O_{\infty}(k^{-\infty})$. Clearly if
$T_k$ is $O_{\infty}(k^{-\infty})$, the compacity of $M_1 \times M_2$
implies that the uniform norm of $\op{T}_k$ is $O(k^{-\infty})$. If $\Pi^1_k \op{T}_k \Pi^2_k = \op{T}_k$ for every $k$, then the converse is true. 

\subsubsection{Fourier integral operators}
If $\om_2$ is the symplectic form of $M_2$, we denote by $M_2^{-}$ the
manifold $M_2$ endowed with the symplectic form $-\om_2$. 
The data to define a Fourier integral operator are a Lagrangian submanifold $\Gamma$ of $M_1 \times M_2^{-}$, a flat unitary section $t_{\Gamma}^{\ka}$ of $L^{\ka}_1 \boxtimes L_2^{-\ka} \rightarrow \Gamma$ and a formal series $\sum \hb^l g_l$ of $\Ci(\Gamma)[[\hb]]$.

By definition $\op{T}_k$ is a {\em Fourier integral operator} associated to $(\Gamma,t_\Gamma^\ka)$ with {\em total symbol} $\sum \hb^l g_l$ if on every compact set $K \subset M_1 \times M_2$ such that $K \cap \Gamma = \emptyset$, 
$$  T_k(x_1,x_2) = O_\infty(k^{-\infty})$$
Furthermore on a neighborhood $U$ of $\Gamma$, 
\begin{gather} \label{kerFIO}
 T_k(x_1,x_2) =  \Bigl( \frac{k}{2\pi} \Bigr)^{n(\Gamma)} E_\Gamma^k(x_1,x_2) f(x_1,x_2,k) + O_{\infty}(k^{-\infty}) \end{gather} 
where 
\begin{enumerate} 
\renewcommand{\labelenumi}{\roman{enumi}}
 \item.   \label{gr}
$E_\Gamma^{\ka}$ is a section of $L^{\ka}_1 \boxtimes L_2^{-\ka} \rightarrow U$ such that $E_\Gamma^{\ka}= t_\Gamma^{\ka}$ over $\Gamma$, and for every holomorphic vector field $Z_1$ of $M_1$ and $Z_2$ of $M_2$ 
\begin{gather*} 
\nabla_{(Z_1,0)} E_\Gamma^{\ka} \equiv 0 \text{ and } \nabla_{(0,\bar{Z}_2)} E_\Gamma^{\ka} \equiv 0 \end{gather*} 
modulo a section which vanishes to any order along $\Gamma$. Furthermore $$|E_\Gamma^{\ka} (x_1, x_2)|< 1$$ if $( x_1, x_2) \notin \Gamma$.  
 \item.  \label{sy}
$(f(.,k ))_k$ is a symbol of $S(U)$ with an asymptotic expansion
$\textstyle{\sum} k^{-l} f_{l}$ such that $$f_l = g_l \text{ over }
\Gamma$$ and  $(Z_1,0). f_{l} \equiv 0 $ and $(0,\bar{Z}_2). f_l
\equiv 0$ modulo a function which vanishes to any order along $\Gamma$ for every holomorphic vector fields $Z_1$ of $M_1$ and $Z_2$ of $M_2$.
\end{enumerate}
$n(\Ga)$ is a real number. Denote by ${\mathcal{F}}(\Gamma, t^{\ka}_\Gamma)$ the set of Fourier integral operators associated to $(\Gamma, t^{\ka}_\Gamma)$. 
\begin{theo} \label{symbFIO}
The map
$ {\mathcal{F}}(\Gamma, t^{\ka}_\Gamma) \rightarrow \Ci(\Gamma)[[\hb]]$ 
which sends an operator into its total symbol is well-defined and onto. Its kernel consists of the operators $O(k^{-\infty})$. 
\end{theo}

The {\em principal symbol} of $\op{T}_k \in {\mathcal{F}}(\Gamma, t^{\ka}_\Gamma)$ is the first coefficient $g_0 \in \Ci(\Gamma)$ of the total symbol. If it doesn't vanish, $\op{T}_k$ is said {\em elliptic}.
In \cite{oim2}, we proved the basic results regarding the composition
properties of this type of Fourier integral operators. 

\subsection{Toeplitz operators} \label{ToNoSc}

The first example of Fourier integral operators are the Toeplitz
operators. 
The diagonal $\De_r$ is a Lagrangian submanifold of $M_r \times M_r^{-}$. Denote by $t^\ka_{\De_r}$ the flat section of $L_r^\kappa \boxtimes L_r^{-\kappa} \rightarrow \Delta_{r}$ such that 
$$t^\ka_{\De_r}(x,x) = z \otimes z^{-1}  \text{ if } z \in L_x^\ka \text{ and }z\neq 0 .$$ 
\begin{defin} ${\mathcal{F}}_r$ is the space of Fourier integral
  operators associated to $(\De_{r},t^\ka_{\De_r})$ with $n(\Delta_r)
  = n_r$ the complex dimension of $M_r$. 
\end{defin} 
By identifying $\De_{r}$ with $M_r$, we consider the total symbols of these operators as formal series of $\Ci(M_r)[[\hb]]$. 
Our main result in \cite{oim1} was the following theorem.
\begin{theo} \label{toepN}
Every Toeplitz operator $(\op{T}_k)$ of $M_r$ is a Fourier integral operator associated to $(\De_{r},t^\ka_{\De_r})$ and conversely. Furthermore, there exists an equivalence of star-products $$E : \Ci(M_r)[[\hb]] \rightarrow \Ci(M_r)[[\hb]]$$ such that if $(\op{T}_k)$ is a the Toeplitz operator with contravariant symbol $\sum \hb^{l} f_l$, then the total symbol of $(\op{T}_k)$ as a Fourier integral operator is $E( \sum \hb^l f_l)$. 
\end{theo} 
The same result holds for the Toeplitz operators of $M$. In the following we use the notations $\De$, $t_{\De}$ and ${\mathcal{F}}$ corresponding to $\De_r$, $t^\ka_{\De_r}$ and ${\mathcal{F}}_r$ on $M$.

\subsection{The $\la$-Toeplitz operators.} \label{laFo}
Recall some notations of section \ref{s2}. The action of $\te \in \T^d$
on $M$ (resp. $L^k$) is denoted by $l_\te$
(resp. ${\mathcal{L}}_\te$). $P$ is the level set $\mu^{-1} (\la)$ and
$p$ is the projection $P \rightarrow M_r$.

Let $\Lambda$ be the moment Lagrangian 
$$\Lambda = \{ (l_{\te}.x,x); \quad x \in P \text{ and } \te \in \T^d \}$$
introduced by Weinstein in \cite{We}. 
$\Lambda$ is a Lagrangian manifold of $M \times M^{-}$. 
Let $t^{\kappa}_\Lambda$ be the section of $L^{\ka} \boxtimes L^{-\ka} \rightarrow \Lambda$ such that  
\begin{gather*} 
 t^{\kappa}_\Lambda ( l_{\te}.x,x ) = {\mathcal{L}}_\te .z \otimes z^{-1} \quad \text{if $z \in L_x^{\ka}$ and $z \neq 0$} \end{gather*}
This is a flat section with constant norm equal to $1$. 
\begin{defin}
${\mathcal{F}}_\la$ is the set of Fourier integral operators $\op{T}_k$ associated to $\Lambda$ and $t^{\kappa}_\Lambda$ with $n(\Lambda) = n - \frac{d}{2}$ and such that 
$$ {\mathcal{L}}_\te^* \op{T}_k = \op{T}_k{\mathcal{L}}_\te^* = \op{T}_k, \quad \forall \te $$
or equivalently $ \Pi_{\la,k} \op{T}_k  \Pi_{\la,k} = \op{T}_k$. 
\end{defin}
We will deduce the following result from the fact that the
algebra ${\mathcal{T}}_\la$ is isomorphic to the algebra
${\mathcal{T}}_r$ of Toeplitz operators of $M_r$ (theorems
\ref{keytheo} and \ref{tptp}). 

\begin{theo} \label{FIO1}
${\mathcal{F}}_\la$ is the algebra ${\mathcal{T}}_\la$ of $\la$-Toeplitz operators.
\end{theo} 

A similar characterization was given by Guillemin-Sternberg for
pseudodifferential operators in \cite{GuSt3}. Their proof starts from the fact that $\Pi_{\la,k}$ is a
Fourier integral operator of ${\mathcal{F}}_\la$. If  $\op{T}_k$ is a
Toeplitz operator of $M$, then the symbolic calculus of Fourier integral operators implies that $\Pi_{\la,k}\op{T}_k \Pi_{\la,k}$ belongs to ${\mathcal{F}}_\la$. This follows essentially from the composition of canonical relations  
$$\Lambda \circ \De \circ \Lambda  = \Lambda .$$ 

In the same way, we can show that ${\mathcal{F}}_\la$ is a $*$-algebra (cf. theorem \ref{theolaT}), define the principal symbol and prove theorem \ref{LaToEq}. The corresponding compositions of canonical relations are
$$ \Lambda \circ \Lambda  = \Lambda, \quad \Lambda^t =\Lambda .$$

The difficulty of this approach is that it uses the properties of composition of the Fourier integral operators. In addition, the composition of $\Lambda$ with itself is not transverse.  But it has the advantage to be more general and can be transposed in other contexts. 

\subsection{The Guillemin-Sternberg isomorphism and its unitarization}

Let $\Theta$ be the Lagrangian submanifold of $M_r \times M^{-}$
$$ \Theta = \{ (p(x),x); \quad x\in P \} $$
Recall that $L_{r}^{\ka} \rightarrow M_r$ is the quotient of $L^\ka \rightarrow P$ by the action of $\T^d$. If $z \in L^{\ka}_x$ where $x \in P$, we denote by $[z] \in L^{\kappa}_{r,p(x)}$ its equivalence class. Then define the section $t^\kappa_\Theta$ of $L^\ka_r \boxtimes L^{-\ka} \rightarrow \Theta$ by 
$$ t^\kappa_{\Theta} (p(x),x) = [z] \otimes z^{-1} \quad \text{if $z \in L_x^{\ka}$ and $z \neq 0$}. $$
It is flat with constant norm equal to $1$. 
\begin{defin} 
  $\flr$ is the space of Fourier integral operators $\op{T}_k$
  associated to $(\Theta, t^\ka_{\Theta})$ with $n(\Theta) = n
  -\frac{3}{4} d$ and such that 
$$ \op{T}_k {\mathcal{L}}_\te^* = \op{T}_k, \quad \forall \te $$
or equivalently $\op{T}_k  \Pi_{\la,k} = \op{T}_k$. 
\end{defin}

\begin{theo} \label{FIO2}
The Guillemin-Sternberg isomorphism $\op{V}_k$ and the unitary operator $\op{U}_k = \op{V}_k  ( \op{V}_k^* \op{V}_k)^{-\frac{1}{2}}$ are elliptic operators of $\flr$.
\end{theo} 

This result is coherent with the fact that $\op{U}_k$ induces an
isomorphism between the algebra of the $\la$-Toeplitz operators and
the algebra of the Toeplitz operators of $M_r$. Indeed, observe that 
$$ \Theta \circ \Lambda \circ \Theta^t = \De_r, \quad \Theta^t \circ \De_r \circ \Theta = \Lambda $$
which corresponds to the equalities
$$ \op{U}_k \op{T}_k \op{U}_k^* = \op{S}_k, \quad\op{U}_k^*  \op{S}_k\op{U}_k= \op{T}_k $$
where $\op{T}_k$ is a $\la$-Toeplitz operator and $ \op{S}_k$ the reduced Toeplitz operator.

To prove theorems  \ref{FIO1} and \ref{FIO2}, we first explain how the spaces $\flr$, $\fl$, $\fr$ of Fourier integral operators are related by the Guillemin-Sternberg isomorphism.  Then we deduce theorems  \ref{FIO1} and \ref{FIO2} from theorem \ref{keytheo}.

\subsection{The relations between $\fl$, $\flr$ and $\fr$} \label{FIOGS}

Since we consider Fourier integral operators which are equivariant, we
need an equivariant version of theorem \ref{symbFIO}. Let us consider
the same data as in section \ref{defFIO}. Let $G$ be a compact Lie
group which acts on $M_1 \times M_2$ preserving the K{\"a}hler structure
of $M_1 \times M_2^{-}$. Assume that this action lifts to
$L_1^\ka \boxtimes L_2^{-\ka}$, preserving the Hermitian structure and
connection.

Suppose that the Lagrangian manifold $\Gamma$ and the section
$t_\Gamma^\ka$ is $G$-invariant. Let us denote by ${\mathcal{F}}_G (
\Gamma, t^{\ka}_\Gamma)$ the space of Fourier integral operators associated to $(\Gamma, t^{\ka}_\Gamma)$ whose kernel is $G$-invariant.

\begin{theo} \label{symbeqFIO}
The total symbol of an operator $\op{T}_k \in  {\mathcal{F}}_G ( \Gamma, t^{\ka}_\Gamma)$ is $G$-invariant. Furthermore the total symbol map $ {\mathcal{F}}_G(\Gamma, t^{\ka}_\Gamma) \rightarrow \Ci_G(\Gamma)[[\hb]]$
is onto. 
\end{theo}

\begin{proof} 
Let $\op{T}_k \in  {\mathcal{F}}_G ( \Gamma, t^{\ka}_\Gamma)$. Its kernel is of the form (\ref{kerFIO}) on a neighborhood of $\Ga$. Hence, 
$$T_k (x_1,x_2) = \Bigl( \frac{k}{2\pi} \Bigr)^{n(\Ga)} f(x_1,x_2,k) t^k_{\Ga} (x_1,x_2) + O(k^{-\infty})$$ 
for every $(x_1,x_2) \in \Ga$. It follows that the total symbol is $G$-invariant. 
Conversely, if the total symbol is $G$-invariant, we can define a $G$-invariant kernel of the form  (\ref{kerFIO}) with a $G$-invariant neighborhood $U$, a $G$-invariant section $E_\Gamma^\ka$ and a $G$-invariant sequence $f(.,k)$. The operator obtained $\op{T}'_k$ does not necessarily satisfy $$\Pi^1_k \op{T}_k' \Pi^2_k =  \op{T}_k'.$$
So we set 
$$ \op{T}_k = \Pi^1_k \op{T}_k' \Pi^2_k.$$ Using that $\Pi^1_k$ and
$\Pi^2_k$ are Fourier integral operators associated to the diagonal of
$M_1$ and $M_2$ respectively, we prove that the kernels of $\op{T}_k$
and $\op{T}_k'$ are the same modulo
$O_\infty(k^{-\infty})$. Consequently $\op{T}_k$ belongs to
${\mathcal{F}}_G ( \Gamma, t^{\ka}_\Gamma)$ and has the required symbol.
\end{proof} 

\begin{cor} \label{ccor}
There is a natural identification between total symbols of
  operators of $\fl$ (resp. ${\mathcal{F}}_{\la,r}$) and formal series
  of $\Ci(M_r)[[\hb]]$.
\end{cor}

\begin{proof} 
By theorem \ref{symbeqFIO}, the total
symbols of the operators of ${\mathcal{F}}_\la$ are the formal series of
$\Ci(\Lambda)[[\hb]]$ invariant with respect to the action of $\T^d
\times \T^d$ on $\Lambda \subset M^2$. The map 
$$ \Lambda \rightarrow M_r, \quad (y,x) \rightarrow p(x) $$ 
is a $(\T^d \times \T^d)$-principal bundle. So there is a one-to-one
correspondence between $\Ci_{\T^d
\times \T^d}(\Lambda)$ and $\Ci(M_r)$. The proof is the same for the
total symbols of the operators of ${\mathcal{F}}_{\la,r}$.
\end{proof}

\begin{theo} \label{keytheo2}
The following maps 
\begin{alignat*}{2} 
 \fl & \rightarrow \flr, \quad  & \op{T}_k &\rightarrow \op{V}_k\op{T}_k \\
 \fl & \rightarrow \fr, \quad  &\op{T}_k &\rightarrow \op{V}_k \op{T}_k \op{V}^*_k
\end{alignat*} 
are well-defined and bijective. Furthermore, if $\op{T}_k \in \fl$,
the total symbols of $\op{V}_k\op{T}_k$, $\op{V}_k \op{T}_k
\op{V}^*_k$ and $\op{T}_k$ are the same with the identifications of
corollary \ref{ccor}.
\end{theo}

\begin{proof} 
These properties follows immediately from the definition of the
 Fourier integral operators. Indeed consider two operators 
$$ \op{T}_k : {\mathcal{H}}_{\la,k} \rightarrow {\mathcal{H}}_{\la,k},
\qquad \op{S}_k : {\mathcal{H}}_{\la,k} \rightarrow
{\mathcal{H}}_{r,k}. $$ 
Extend them to the space of $L^2$ sections in such a way that they vanish on the orthogonal of ${\mathcal{H}}_{\la,k}$. Then $\op{S}_k = \op{V}_k  \op{T}_k$ or equivalently $\op{T}_k =\op{W}_k \op{S}_k$ if and only if the kernels of $\op{T}_k$ and $\op{S}_k$ satisfy\begin{gather} \label{pass}
 T_k = \Bigl( \frac{k}{2 \pi} \Bigr)^{\frac{d}{4}} ( p_\C \boxtimes \op{Id} )^* S_k \quad \text{ over } P_\C \times M. \end{gather}
This follows from the definition of $V_k$ (cf. definition \ref{defGSmap}). 
Furthermore by proposition \ref{pcompl}, the kernel of $T_k$ vanishes over $P_{\C}^c \times M$. 
So we may recover the kernel of $\op{T}_k$ from the kernel of
$\op{S}_k$ and conversely. Using this we can directly check that $\op{T}_k \in \fl$ if and only if $\op{S}_k \in \flr$.

To do this observe that the data which define the operators of $\fl$ and $\flr$ are related in the following way:
$$ (p \times \op{Id} )^{-1}( \Theta) =\Lambda, \qquad   (p \boxtimes \op{Id})^* t_\Theta^\kappa = t_\Lambda^\kappa $$
Using that $p_\C :P_\C \rightarrow M_r$ is a holomorphic map, it comes
that 
$$ (p_\C \boxtimes \op{Id} )^*  E_{\Theta}^\kappa =
E_{\Lambda}^\kappa .$$
That $\op{T}_k \rightarrow \op{V}_k\op{T}_k$ is well defined, bijective and preserves the total symbols follows easily.

For the second map, we can proceed in a similar way. Let $\op{S}_k$ be
such that $\Pi_{r,k} \op{S}_k \Pi_{r,k}= \op{S}_k$. Then computing
successively  the Schwartz kernels of $\op{W}_k \op{S}_k$,
$(\op{W}_k\op{S}_k)^*$, $\op{W}_k (\op{W}_k\op{S}_k)^*$ and using
$$( \op{W}_k (\op{W}_k\op{S}_k)^*)^* =  \op{W}_k\op{S_k} \op{W}_k^*,$$ we obtain that the kernels of $\op{T}_k =  \op{W}_k\op{S_k} \op{W}_k^*$ and $\op{S}_k$ satisfy  
$$ T_k(x,y) = \Bigl( \frac{k}{2 \pi} \Bigr)^{\frac{d}{2}} e^{-k \ph (y)} ( p_\C \boxtimes p_\C )^* S_k (x,y) \quad \text{ over } P_\C \times P_\C.$$
where the function $\ph$ has been defined in proposition \ref{pcontrolnorm}. 
\end{proof}

\subsection{Proofs of theorems \ref{FIO1} and \ref{FIO2}}

Recall that by theorem \ref{keytheo}, there is a bijection from ${\mathcal{T}}_\la$ onto ${\mathcal{T}}_r$
\begin{gather} \label{asd1}
 \laT \rightarrow \rT, \quad \op{T}_k \rightarrow \op{W}_k^* \op{T}_k \op{W}_k. \end{gather}
Furthermore, we know that ${\mathcal{T}}_r = {\mathcal{F}}_r$ by
theorem \ref{toepN}. In theorem \ref{keytheo2}, we proved that the map
\begin{gather} \label{asd2}
 \fl \rightarrow \fr,   \quad \op{T}_k \rightarrow \op{V}_k \op{T}_k \op{V}_k^* \end{gather}
is a bijection. Let us deduce that $ {\mathcal{T}}_\la = {\mathcal{F}}_\la$. 

Let $\op{T}_k$ belong to $\fl$. By (\ref{asd2}), $\op{V}_k\op{T}_k \op{V}_k^*$ belongs to ${\mathcal{T}}_r$. By (\ref{asd1}), $\op{W}_k^* \op{W}_k$ belongs to ${\mathcal{T}}_r$. Since the product of Toeplitz operators is a Toeplitz operator, $$ (\op{W}_k^* \op{W}_k)(\op{V}_k\op{T}_k \op{V}_k^*) (\op{W}_k^* \op{W}_k)= \op{W}_k^* \op{T}_k\op{W}_k$$
 belongs to ${\mathcal{T}}_r$. By (\ref{asd1}) $\op{T}_k$ belongs to ${\mathcal{T}}_\la$.

Conversely assume that $\op{T}_k$ belongs to ${\mathcal{T}}_\la$. By (\ref{asd1}), $\op{W}_k^* \op{T}_k\op{W}_k$  belongs to ${\mathcal{T}}_r$. Write 
\begin{align*} 
 \op{V}_k \op{T}_k \op{V}_k^* =  & ( \op{V}_k \op{V}_k^* ) (\op{W}_k^* \op{T}_k \op{W}_k) ( \op{V_k}\op{V}_k^*) \end{align*}
Observe that $\op{V}_k \op{V}_k^*$ is the inverse of $\op{W}_k^* \op{W}_k$ in the sense of Toeplitz operators, i.e.   
\begin{gather*} 
\Pi_{r,k} (\op{V}_k \op{V}_k^*)  \Pi_{r,k} = \op{V}_k \op{V}_k^* , \\
  (\op{V}_k \op{V}_k^*)  (\op{W}_k^* \op{W}_k)= (\op{W}_k^* \op{W}_k) (\op{V}_k \op{V}_k^*)=  \Pi_{r,k}. \end{gather*} 
Since $(\op{W}_k^* \op{W}_k)$ is a Toeplitz operator with a non-vanishing symbol by theorem \ref{keytheo}, $\op{V}_k \op{V}_k^*$ is a Toeplitz operator. Consequently $ \op{V}_k \op{T}_k \op{V}_k^*$ belongs to ${\mathcal{T}}_r$ and by (\ref{asd2}), $\op{T}_k$ belongs to $\fl$.

Let us prove theorem \ref{FIO2}. By theorem \ref{FIO1}, $\Pi_{\la,k}$ belongs to $\fl$. So theorem \ref{keytheo2} implies that $\op{V}_k = \op{V}_k \Pi_{\la,k}$ belongs to $\flr$. 
Let us consider now $\op{U}_k$. We have 
$$ \op{U}_k = (\op{W}_k^* \op{W}_k )^{-\frac{1}{2}} \op{W}_k^* = \op{V}_k (  \op{W}_k  (\op{W}_k^* \op{W}_k )^{-\frac{1}{2}} \op{W}_k^*)$$ 
As we saw in the proof of theorem \ref{tptp}, $ (\op{W}_k^* \op{W}_k )^{-\frac{1}{2}}$ belongs to ${\mathcal{T}}_r$. So by theorem \ref{keytheo},  $\op{W}_k  (\op{W}_k^* \op{W}_k )^{-\frac{1}{2}} \op{W}_k^*$ belongs to ${\mathcal{T}}_\la = {\mathcal{F}}_\la$. And theorem \ref{keytheo2} implies that $\op{U}_k$ belongs to $\flr$.


\section{Toeplitz operators on orbifold} \label{ToepOrb}

In this part, we prove the basic results about the Toeplitz operators
on the orbifold $M_r$. We describe their kernels as Fourier integral
operators associated to the diagonal, prove that the set of Toeplitz
operators is an algebra and describe the associated symbolic
calculus. Finally we compute the asymptotic of the density of states of a Toeplitz operator.

\subsection{Schwartz kernel on orbifold} 
Let us introduce some notations and state some basic facts
about kernels of operators on orbifolds. Let $X$ be
a reduced orbifold with a vector bundle $E \rightarrow X$. So every
chart $(|U|,U, G,\pi_U)$ of $X$  is endowed with a $G$-bundle $E_U \rightarrow
U$. $G$ acts effectively on $U$. If $g \in G$, we denote by $a_g: U \rightarrow U$ its action on
$U$ and by ${\mathcal{A}}_g : E_U \rightarrow E_U$ its lift.
If $s$ is a section of $E \rightarrow X$, we denote by $s_U$ the
corresponding invariant section of $E_U$. 

As in the manifold case, we can define the dual bundle of $E$, the tensor
product of two bundles over $X$, the orbifold $X^2$ and the
bundle $E \boxtimes E^* \rightarrow X^2$. 
Let $\delta$ be a volume form of $X$. Then every section $T$ of $E \boxtimes
E^*$ defines an operator $\op{T}$ in the following way. Consider two charts
$(|U|,U, G,\pi_U)$ and $(|V|,V, H ,\pi_V)$ of $X$. Then $T$ is given
over $|U| \times |V|$ by a $(G \times H)$-invariant section $T_{UV}$
of $E_U \boxtimes E_V^*$. If $s$ is a section of $E$ with compact support
in $|V|$, then 
$$ (\op{T} s )_U = \frac{1}{\# H} \int_V T_{UV}. s_V \; \delta_V. $$
Assume that $E$ is Hermitian and define the scalar product of
sections of $E$ by using $\delta$. If $\op{T}$ is an operator which acts on $\Ci( X,
E)$, vanishing over the orthogonal of a finite dimensional subspace of $\Ci( X,
E)$, it is easily proved that $\op{T}$ has a Schwartz
kernel $T$. It is unique.
Furthermore, if $E$ has rank one, the trace of $\op{T}$ is given by 
$$ \op{Tr} ( \op{T}) = \int_X \Delta^* T \; \delta $$
where $\Delta : X \rightarrow X^2$ is the diagonal map. Since $\Delta$ is a
good map
(in the sense of \cite{ChRu}), the
pull-back $\Delta^* (E \boxtimes E^*)$ is well-defined. It is
naturally isomorphic to $ E \otimes E^* \rightarrow X$. So $\Delta^*
T$ is a section of $E \otimes E^* \simeq \C$. Finally the previous
integral is defined with the orbifold convention: if $
(|U|,U, G,\pi_U)$ is a chart of $X$ and $\Delta^* T$ has support in $|U|$, it is given by 
$$ \frac{1}{\# G}  \int_U (\Delta^* T)_U \; \delta_U.$$

Note also that it is false that every section of $E \otimes E^*$
is the pull-back by $\Delta$ of a section of $E \boxtimes E^*$. For instance when $E$
is a line bundle, $E\otimes E^* \simeq \C$ has nowhere
vanishing sections, whereas it may happen that every section of $E
\boxtimes E^*$ vanishes at some point $(x,x)$ of the diagonal. This
explains some complications in the description of the kernels of
Toeplitz operators. 

\subsection{The algebra $\fr$} \label{secork}

Recall that $(M_r, \om_r)$ is a compact K{\"a}hler reduced orbifold with a
prequantum bundle $L^\ka_r \rightarrow M_r$ whose curvature is $- i \ka 
\om_r$. Let us consider a family $(\op{T}_k)_k$ of operators, with Schwartz
kernels 
$$T_k \in \Ci(M_r^2, L_{r}^k \boxtimes L^{-k}_r), \quad k = \ka, 2\ka,
3 \ka,...$$
As in the manifold case (cf. section \ref{ToNoSc}), the definition of the operators of $\fr$ consists in two parts. The first assumption is 
\begin{AS} \label{AsTo1} 
$T_k \text{ is } O_{\infty}(k^{-\infty})\text{ on every compact set }K
\subset M_r^2 \text{ such that } K \cap \De_r = \emptyset .$
\end{AS}
The description
of $T_k$ on a neighborhood of the diagonal doesn't generalize
directly from the manifold case, because the definition of the section $E^\ka_{\De_r}$
doesn't make sense. Fortunately we can keep the same ansatz on the
orbifold charts. If $(|U|,U, G,\pi_U)$ is an orbifold chart of
$M_r$, we assume that:

\begin{AS} \label{AsTo2} 
There exists a section $T'_{k,U}$ of $L^k_{r,U} \boxtimes L^{-k}_{r,U}$ invariant with respect to the diagonal action of $G$ and 
of the form 
\begin{gather} \label{kerFIO2} 
 T'_{k,U}(x,y) = \Bigl( \frac{k}{2\pi} \Bigr)^{n_r} E^k_{\De_U}(x,y) f(x,y,k) + O_\infty(k^{-\infty}) \end{gather}
on a neighborhood of the diagonal $\De_U$, where $E^\ka_{\De_U}$ and
$f(.,k)$ satisfy the assumptions
(\ref{kerFIO}.  {\romannumeral 1})
and (\ref{kerFIO}. {\romannumeral 2}) 
with $(\Ga, t_\Ga^\ka) = (\De_{U}, t_{\De_{U}}^\ka)$, such that  
\begin{gather} \label{Tmoy}
 T_{k,UU} = \sum_{g \in G}  ({\mathcal{A}}_g \boxtimes \op{Id} )^* T'_{k,U} .\end{gather}
\end{AS}
Observe that we can use $T'_{k,U}$ instead of $T_{k,UU}$ to compute  $(\op{T}_k \Psi)_U$ when $\Psi$ has compact support in $U$:
\begin{xalignat}{2} \notag   
  (\op{T}_k \Psi)_U (x) = & \frac{1}{\# G} \int_U    T_{k,UU} (x,y) \Psi_U(y)
  \delta_{U} (y) .
\\ \label{locac}
 = & \int_U    T'_{k,U} (x,y) \Psi_U(y)
  \delta_{U} (y) .
\end{xalignat}
This follows from (\ref{Tmoy}) and the fact that $\Psi_U$ is $G$-invariant and $\op{T}'_{k,U}$ is invariant  with respect to the diagonal action.

\begin{defin}
$\fr$ is the set of operators $(\op{T}_k)$ such that $$\Pi_{r,k}  \op{T}_k \Pi_{r,k} = \op{T}_k$$ and whose Schwartz kernel satisfies assumptions \ref{AsTo1} and \ref{AsTo2} for every orbifold chart $(|U|,U, G,\pi_U)$. 
\end{defin}

The basic result is the generalization of the theorem of Boutet de Monvel and Sj{\"o}strand on the Szeg{\"o} projector. 
\begin{theo} \label{BSor}
The projector $\Pi_k$ is an operator of ${\mathcal{F}}_r$. \end{theo}
This theorem will be proved in section \ref{prep}. Let us deduce from it the properties of ${\mathcal{F}}_r$. 
First we define the total symbol map 
$$ \si: \fr \rightarrow \Ci(M_r) [[\hb]], $$
which sends $\op{T}_k$ into the formal series $\sum \hb^l g_l$ such that 
\begin{gather} \label{dsp}
 T_k (x,x) = \bigl(\tfrac{k}{2\pi} \bigr)^{n_r} \textstyle{\sum} k^{-l} g_l(x) + O(k^{-\infty}) \end{gather}
for every $x$ which belongs to the principal stratum of $M_r$.

\begin{prop} \label{totsymb}
$\si$ is well-defined and onto. Furthermore $\si (\op{T}_k) =0$ if and
only if $\op{T}_k$ is smoothing in the sense of section \ref{smooth}, i.e. its
kernel is  $O_\infty(k^{-\infty})$.
\end{prop}
\begin{proof} First let us prove that $\si$ is well-defined.   By (\ref{Tmoy}),
$$  T_{k,UU}(x,x) =  \textstyle{\sum}_{g \in G} ({\mathcal{A}}_g \otimes \op{Id} )^* T'_{k,U}(a_g.x,x) $$
Assume that $\pi_U(x)$ belongs to the principal stratum of $M_r$. So if $g\neq \op{id}_G$, $a_g.x \neq x $ which implies that $T'_k(a_g.x,x) = O(k^{-\infty})$. Consequently (\ref{kerFIO2}) leads to
$$  T_{k,UU}(x,x) = \bigl(\tfrac{k}{2\pi} \bigr)^{n_r} f(x,x,k) + O(k^{-\infty}).$$
 This proves the existence of the asymptotic expansion (\ref{dsp}) and
 that the $g_l$ extend to $\Ci$ functions on $M_r$. Furthermore, since
 the principal stratum is dense in $M_r$, these functions are uniquely
 determined by the kernel $T_k$. If $\op{T}_k$ is smoothing, they vanish.
Conversely, if the functions $g_l$ vanish, $T'_{k,U}$ is $O_\infty(k^{-\infty})$ and the same holds for $T_{k,UU}$. So the kernel of $\si$ consists of the $\op{T}_k$ such that $T_k$ is   $O_\infty(k^{-\infty})$. 

Let us prove that $\si$ is onto. Let $\sum \hb^l g_l$ be a formal
series of $\Ci(M_r)[[\hb]]$. First we construct a Schwartz kernel
$T_{k}$ satisfying assumptions \ref{AsTo1}, \ref{AsTo2} and (\ref{dsp}). To do this, we introduce on every orbifold chart a kernel $T_{k,U}'$ of the form (\ref{kerFIO2}) where the functions $f(.,k)$ are such that 
$$ f(.,k) = \textstyle{ \sum} k^{-l} g_{l,U} + O(k^{-\infty})$$ 
with the $g_{l,U}\in \Ci(U) $ corresponding to the $g_l$. The existence of $T_{k,U}'$ is a consequence of the Borel lemma as in the manifold case. Furthermore since the functions $g_{l,U}$ are $G$-invariant and the diagonal action of $G$ preserves $\De_U$ and $t^\ka_{\De_{U}}$, we can obtain a $G$-invariant section $T_{k,U}'$. 

Then we piece together the $T_{k,UU}$ by using a partition of unity
subordinate to a cover of $M_r$ by orbifold charts. To do this, we
have to check that two kernels $T_{k,U_1 U_1}$ and $T_{k,U_2U_2}$ on two
orbifold charts 
$$(|U_1|,U_1, G_1,\pi_{U_1}), \quad (|U_2|,U_2,G_2,\pi_{U_2})$$
 define the same section over $|U_1| \cap |U_2|$ modulo $O_{\infty}(k^{-\infty})$. Recall that the compatibility between orbifold charts is expressed by using orbifold charts $(|U|,U, G,\pi_U)$ which embed into $(|U_i|,U_i, G_i,\pi_{U_i})$. Denote by $\rho_i: G \rightarrow G_i$ the injective group homomorphism and by $j_i : U \rightarrow U_i$ the $\rho_i$-equivariant embeddings. We have 
\begin{gather} \label{ghj}
 \forall \; g \in G_i, \; \forall x \in j_i(U), \; a_g.x \in   j_i(U) \Rightarrow g \in  \rho_i(G). \end{gather}
We have to prove that 
\begin{gather} \label{ckl}
 (j_i \boxtimes j_i )^* T_{k,U_iU_i} = T_{k,UU} + O_{\infty}(k^{-\infty}) \end{gather}
By (\ref{ghj}), if $g \in G_i - \rho_i(G)$, then $a_g ( j_i(U)) \cap   j_i(U) = \emptyset$. Since $T'_{k,U_i}$ is $O_\infty(k^{-\infty})$ outside the diagonal, this implies 
$$  T_{k,U_iU_i}\eval{j_i(U) \times j_i(U)} = \sum_{g \in \rho_i(G) }  ({\mathcal{A}}_g \boxtimes \op{Id} )^* T'_{k,U_i}\eval{j_i(U) \times j_i(U)} + O_{\infty}(k^{-\infty}) $$
Furthermore, since $j_i^* g_{l,U_i} = g_{l,U}$, we have 
$$ (j_i \boxtimes j_i)^* T'_{k,U_i} = T'_{k,U} +  O_{\infty}(k^{-\infty}) $$
Both of the previous equation lead to (\ref{ckl}) by using that the map $j_i$ are $\rho_i$-equi\-va\-riant.

The section $T_{k}$ obtained is the Schwartz kernel of an operator $\op{T}_k$. We don't have necessarily $\Pi_k \op{T}_k \Pi_k = \op{T}_k$, but only
\begin{gather} \label{fr}
 \Pi_k \op{T}_k \Pi_k \equiv \op{T}_k  \end{gather}
modulo an operator whose kernel is $O_{\infty}(k^{-\infty})$. So we
replace $\op{T}_k$ with $\Pi_k \op{T}_k \Pi_k$. The proof of
(\ref{fr}) is a consequence of theorem \ref{BSor} and is similar to
the manifold case. Actually, if two operators $\op{R}_k$ and $
\op{S}_k$  satisfy assumptions \ref{AsTo1} and \ref{AsTo2},  the kernel $T_k$ of their product is given on a orbifold chart by 
$$ T_{k,UU} = \sum_{g \in G}  ({\mathcal{A}}_g \boxtimes \op{Id} )^* T'_{k,U} $$
where 
\begin{gather} \label{locprod}
  T'_{k,U}(x,y) = \int_U S'_{k,U}(x,z). R'_{k,U}(z,y)\ \de_{U}(z) + O_\infty(k^{-\infty}) \end{gather}
This follows from equation (\ref{locac}). 
\end{proof}

\begin{theo} $\fr$ is a $*$-algebra and the induced product of total symbols is a star-product of $\Ci(M_r)[[\hb]]$. Every operator of $\fr$ is of the form
\begin{gather} \label{tM}
 \Pi_{r,k} M_{f(.,k)} \Pi_{r,k} + O(k^{-\infty}) \end{gather} 
where $f(.,k)$ is a symbol of $S(M_r)$ and conversely. The map 
$$ E : \Ci(M_r)[[\hb]] \rightarrow \Ci(M_r)[[\hb]]$$
which sends the formal series $\sum \hb^l f_l$ corresponding to the
multiplicator $f(.,k)$ into the total symbol of the operator
(\ref{tM}), is well defined. It is an equivalence of star-products. 
\end{theo}

\begin{proof} 
As we have seen in the previous proof, the composition of kernels of
operators of $\fr$ corresponds in a orbifold chart to a composition on
a manifold (cf. equation (\ref{locprod})). So the proof in the
manifold case \cite{oim1} extends directly. Let us recall the main
steps. We first prove that $\fr$ is a $*$-algebra and compute the
product of the total symbols by applying the stationary phase lemma to
the composition of the kernels. In the same way, we can compute the
kernel of the operator (\ref{tM}), since $\Pi_{r,k}$ belongs to $\fr$
by theorem \ref{BSor}. As a result of the computation, this operator
belongs to ${\mathcal{F}}_r$ and we obtain that the map $E$ is an
equivalence of star-product. From this we deduce that conversely every operator of $\fr$ is of the form (\ref{tM}).
\end{proof} 

\subsection{Spectral density of a Toeplitz operator}

Let $\op{T}_k$ be a self-adjoint Toeplitz operator over $M_r$. Denote
by $d_k$ the dimension of ${\mathcal{H}}_{r,k}$ and by $E_1 \leqslant
E_2 \leqslant ... \leqslant E_{d_k}$ the eigenvalues of $\op{T}_k$. The spectral density of $\op{T}_k$ is the measure of $\R$
$$   \mu_{\op{T}_k}(E) = \sum_{i=1 }^{d_k} \delta (E - E_i) $$
Let $f$ be a $\Ci$ function on $\R$. We will estimate $\langle \mu_{\op{T}_k} , f \rangle $. The first step is to compute the operator $f(\op{T}_k)$. 

\begin{theo} 
 If $f$ is a  $\Ci$ function on $\R$, then $f(\op{T}_k)$ is a Toeplitz operator. Furthermore if $g_0$ is the principal symbol of $\op{T}_k$ then $f(g_0)$ is the principal symbol of $f(\op{T}_k)$.
\end{theo}
 
The proof is similar to the manifold case (cf. proposition 12 of \cite{oim1}). Now we have 
$$  \langle \mu_{\op{T}_k} , f \rangle =  \sum_{i=1 }^{d_k} f(E_i) =
  \op{Tr} f(\op{T}_k) $$
By the previous theorem, it suffices to estimate the trace of a
  Toeplitz operator. Recall that 
$$ \op{Tr} \op{T}_k = \int_{M_r} T_k (x,x)  \delta_{M_r}(x) $$
Let us begin with the computation in a orbifold chart $(|U|,U,
G,\pi_U)$. Let $\eta$ be a $\Ci$ function of $M_r$ whose support is
included in $|U|$. It follows from assumption \ref{AsTo2} that 
$$ \int_{M_r} \eta(x) T_k (x,x)  \delta_{M_r}(x) =  \frac{1}{\# G} \sum_{g \in G} I(g,U)$$ 
where $$I(g,U) = \int_U  \eta_U(x) ({\mathcal{A}}_g \boxtimes \op{Id} )^* T'_{k,U} (x,x)  \delta_{U}(x) .$$

Let $U^g$ be the fixed point set of $g$,
\begin{gather} \label{defUg} 
U^g= \{ x \in U; \; a_g.x =x \}. \end{gather}
 Assume that $U^g$ is connected. Then $U^g$ is K{\"a}hler submanifold of $U$. Denote by $d(g)$ its complex codimension. 

Let $y \in U^g$. $a_g: U \rightarrow U$ induces a linear transformation on the normal
space $N_y= T_y U/ T_y U^g$. It is a unitary map with eigenvalues 
$$b_1(g),...,b_{d(g)}(g)$$ 
on the unit circle and not equal to $1$. Furthermore, ${\mathcal{A}}_g$ acts on the fiber of $L^\ka_{r,U}$ at $y$ by multiplication by $c^\ka(g)$, where $c(g)$ is on the unit circle. Since $U^g$ is connected, the complex numbers $b_i(g)$ and $c(g)$ do not depend on $y$. 

\begin{lemme}  \label{ll1}
The integral $I(g,U)$ admits an asymptotic expansion
$$ I(g,U) = \Bigl( \frac{k}{2\pi} \Bigr)^{n_r - d(g) }  c(g)^{-k}  \sum_{l=0}^{\infty} k^{-l} I_l(g,U) + O(k^{-\infty}) .$$
 The first coefficient is given by 
$$ I_0(g,U) = \Biggl( \prod_{i=1}^{d(g) } (1 - b_i(g) )^{-1} \Biggr)  \int_{U^g} \eta_{U}(x) f_{0,U} (x) \delta_{U^g}(x) $$ 
where $f_0$ is the principal symbol of $\op{T}_k$ and $\delta_{U^g}$ is the Liouville measure of $U^g$. 
\end{lemme}

\begin{proof} 
Since $|T'_k(x,y)|$ is $O(k^{-\infty})$ if  $x \neq y$, we can restrict the integral over a neighborhood of $U^g$. Let $y \in U^g$. Let $(w^i)$ be a system of complex coordinates of $U$ centered at $y$. Since $g$ is of finite order, we can linearize the action of $a_g$. So we can assume that 
$$ a_g : (w_i)_{i=1,...,n_r} \rightarrow (w_1 b_1, ...,w_{d} b_d, w_{d+1},..., w_{n_r}) $$
To simplify the notations we denoted by $d$, $b_i$ the numbers $d(g)$, $b_i(g)$. In the same way, we can choose a local holomorphic section $s^\ka_r$ of $L^\ka_{r,U}$, which doesn't vanish on a neighborhood of $y$ and such that 
$$ {\mathcal{A}}_g^* s^\ka_r = c^{-\ka} (g) s^\ka_r$$ 
Let $H_r$ be the function such that $ (s^\ka_r , s^\ka_r) = e^{-\ka H_r}$. Observe that $H_r (a_g.x) =H_r(x)$. 
Denote by $t^\ka$ the unitary section $e^{\ka H/2} s^\ka$. 
Introduce a function $\tilde{H}_r(x,y)$ defined on $U^2$ such that $\tilde{H}_r(x,x) = H_r(x)$ and
$$(0,\bar{Z}).\tilde{H}_r \equiv 0 \text{ and } (Z,0).\tilde{H}_r
 \equiv 0$$
modulo $O(|x-y |^{\infty})$ for every holomorphic vector field
$Z$ of $U$.
Using that 
$$  \nabla t^\ka = \frac{\ka}{2} ( \bar{\partial} H - \partial H)
\otimes t^\ka $$
we compute the section $ E^k_{\De_U}$ of \eqref{kerFIO2}
(cf. \cite[proposition 1]{oim1}), and obtain
$$ T'_k(x,y) = \Bigl( \frac{k}{2\pi} \Bigr)^{n_r} e^{- k \bigl( \frac{1}{2}(H_r(x) + H_r(y)) - \tilde{H}_r(x,y) \bigr)} f(x,y,k) t^k_r(x) \otimes t^{-k}_r(y) $$ 
Consequently the integral $I(g,U)$ is equal to 
$$ \Bigl( \frac{k}{2\pi} \Bigr)^{n_r} c^{-k} (g) \int_{U} e^{-k \phi_g (x)} f(a_g.x,x,k) \delta_{U} (x) $$
where $\phi_g (x) = H_r(x) - \tilde{H}_r(a_g.x,x)$. We estimate it by applying the stationary phase lemma. $\phi_g$ vanishes along $U^g$. Using that $a_g^* H_r = H_r$, we obtain for $i,j=1,...,d$, 
$$ \partial_{w_i} \phi_g = \partial_{\bar{w}_i} \phi_g = \partial_{w_i}\partial_{w_j} \phi_g  = \partial_{\bar{w}_i}\partial_{\bar{w}_j} \phi_g = 0$$
along $U^g$ for $i,j=1,...,d$. Denote by $H_{i,j}$ the second derivative $\partial_{w_i}\partial_{\bar{w}_j} H_r$. Then 
$$ \partial_{w_i}\partial_{\bar{w}_j} \phi_g = (1 - b_i ) H_{i,j} $$
along $U^g$ for $i,j=1,...,d$. Furthermore 
$$ \delta_{U} (x) = \op{det} [ H_{i,j} ]_{i,j=1,...,n_r} |dw_1d\bar{w}_1...dw_{n_r}d\bar{w}_{n_r}|$$
which leads to the result. 
\end{proof} 

The next step is to patch together these local contributions. This
involves a family of orbifolds associated to $M_r$ which appears also
in the Riemann-Roch theorem or in the definition of orbifold
cohomology groups (cf. the associated orbifold of \cite{Me}, the
twisted sectors of \cite{ChRu}, the inertia orbifold of
\cite{Mo}). The description of these orbifolds in the general case is
rather complicated. Here, $M_r$ is the quotient of $P$ by a torus
action, which simplifies the exposition (cf. the
appendix of \cite{GuGiKa}).

Consider the following set 
$$ \tilde{P} = \{ (x, g ) \in P \times \T ^d ; \; l_g.x =x \}$$
To each connected component $C$ of $ \tilde{P}$ is associated a
element $g$ of $\T^d $ and a {\em support} $\bar{C} \subset P$ such that $C = \bar{C}
\times \{ g \}$. $\bar{C}$ is a closed submanifold of $P$ invariant
with respect to the action of $\T^d$. The quotient $$ F:=  \bar{C} /
\T^d $$ is a compact orbifold, which embeds into $M_r$. Since $\T^d$
doesn't necessarily act effectively on $\bar{C}$, $F$ is not
in general a reduced orbifold. Denote by $m(F)$ its multiplicity. 

Let $(x,g) \in C$ and denote by $G$ the isotropy group of $x$. Let $U \subset P$
be a slice at $x$ of the $\T^d$-action. Let $|U| = p(U)$ and $\pi_U$
be the projection $U \rightarrow |U|$. Then $(|U|,U, G, \pi_U)$ is an
orbifold chart of $M_r$. Introduce as in (\ref{defUg}) the subset $U^g$ of
$U$ and assume it is connected. Then $U^g = U \cap \bar{C}$. Let
$|U^g| = p(U^g)$ and $\pi_{U^g}$ be the projection $U^g \rightarrow
|U^g|$. Then $(|U^g|,U^g, G, \pi_{U^g})$ is an orbifold chart of
$F$. So $I_0(g,U)$ in lemma \ref{ll1} is given by an integral over
$F$.  
Furthermore, since $F$ is connected, there exists
complex numbers 
$$b_1(g,F),..., b_{d(F)}(g,F) \text{  and } c(g,F)$$ on the unit
circle corresponding to the numbers defined locally, with $d(F)$ the
codimension of $F$ in $M_r$.   Observe also
that $F$ inherits a K{\"a}hler structure. 
 
Denote by ${\mathcal{F}}$ the set 
$$ {\mathcal{F}} := \{ \bar{C}/ \T^d; \; \bar{C} \times \{ g \} \text{ is a component of $\tilde{P}$} \} $$ 
For every $F \in {\mathcal{F}}$, let $\T^d_F$ be the set of $g \in
\T^d$ such that $F = \bar{C} / \T^d$ and $\bar{C} \times \{ g \}$ is
a component of $\tilde{P}$. The point is that two components of
$\tilde{P}$ may have the same support. Since the set of components of
$\tilde{P}$ is finite, the various sets ${\mathcal{F}}$ and $\T^d_F$ are finite.

\begin{theo} \label{Tra}
Let $\op{T}_k$ be a self-adjoint Toeplitz operator on $M_r$ with principal symbol $g_0$. Let $f$ be a $\Ci$ function on $\R$. Then $\langle \mu_{\op{T}_k} , f \rangle$ admits an asymptotic expansion of the form
$$ \sum_{F \in {\mathcal{F}}}   \Bigl( \frac{k}{2\pi} \Bigr)^{n_r - d(F) }  \sum_{g \in \T^d_F}  c(g,F)^{-k}  \sum_{l=0}^{\infty}  k^{-l} I_l(F,g) + O(k^{-\infty})
$$
where the coefficients $I_l(F,g)$ are complex numbers. Furthermore
$$ I_0(F,g) =  \frac{1}{m(F)} \Biggl( \prod_{i=1}^{d(F) } (1 - b_i(g,F) )^{-1} \Biggr) \int_F f(g_0) \delta_{F}$$
where $\delta_{F}$ is the Liouville measure of $F$.
\end{theo} 

\begin{rem}
Each orbifold $F \in {\mathcal{F}}$
is the closure of a strata of $M_r$ (cf. appendix of
\cite{GuGiKa}). For instance, $M_r$ itself belongs to ${\mathcal{F}}$
and is the closure of the principal stratum of $M_r$. Note that
$\T^d_{M_r} = \{ 0 \}$ and that the other suborbifolds of ${\mathcal{F}}$
have a positive codimension. So at first order, 
$$ \langle \mu_{\op{T}_k} , f \rangle =  \Bigl( \frac{k}{2\pi}
\Bigr)^{n_r} \int_{M_r}  f(g_0) \delta_{M_r} + O( k^{n_r-1}).$$
\end{rem}

\begin{rem} \label{remRKK}
Riemann-Roch-Kawasaki theorem gives the dimension of
${\mathcal{H}}_{r,k}$ in terms of characteristic forms, when $k$ is
sufficiently large: 
$$ \op{dim} {\mathcal{H}}_{r,k} = \sum_{F \in {\mathcal{F}}}  \sum_{g
  \in \T^d_F} \frac{1}{m(F)} \int_{F} \frac{ \op{Td}(F) \op{Ch}( L_r^k,F, g) }{
  \op{D}(N_{F}, g)} $$
For the definition of these forms, we refer to theorem 3.3 of
  \cite{Me}. Let us compare this with the estimate of $\op{Tr}
  (\Pi_k)$ given by theorem \ref{Tra}. First,
 $\op{Ch}( L_r^k, F,g)\in \Om(F)$ is a twisted characteristic form associated to
  the pull-back of $L_r^k$ by the embedding $F \rightarrow M_r$. At
  first order
$$ \op{Ch}( L_r^k, F,g) =  c(g,F)^{k}  \Bigl( \frac{k}{2\pi} \Bigr)^{n_F}
\frac{ \om_F^{\wedge n_F}} { n_F !} + O( k ^{n_F -1}) $$
where $n_F =  n_r - d(F)$ and  $\om_F$ is the symplectic form of
$F$. $\op{D}(N_{F}, g)$ is a
twisted characteristic form associated to the normal bundle $N_{F}$ of the
embedding $F \rightarrow M_r$, 
$$ 
\op{D}(N_{F}, g) \equiv \prod_{i=1}^{d(F) } (1 - b_i(g,F)) \mod
\Om^{\bullet \geqslant 2}(F)  $$
$\op{Td}(F)\equiv 1$ modulo $\Om^{\bullet \geqslant 2}(F)$ is the Todd form of
$F$. Hence we recover the leading term $I_0(F,g)$ in theorem
\ref{Tra}. \qed
\end{rem}
 
\begin{rem} (harmonic oscillator). We can describe explicitly  the
  various term of theorem \ref{Tra} and deduce theorem \ref{int2}. 
 Denote by ${\mathcal{P}}$ the set of greatest common divisors of the families $({\mathfrak{p}}_i)_{i \in I}$ where $I$ runs over the subsets of $\{1,...,n\}$. For every ${\mathfrak{p}} \in {\mathcal{P}}$, define the  subset $I({\mathfrak{p}})$ of $\{1,...,n\}$ 
$$ I({\mathfrak{p}}) := \{ i ; \;  \text{${\mathfrak{p}}$ divides ${\mathfrak{p}}_i$ } \} $$
and the symplectic subspace $\C_{\mathfrak{p}}$ of $\C^n$
$$  \C_{\mathfrak{p}} := \{ z \in \C^n; \;  z_i = 0 \text{ if } i
\notin I({\mathfrak{p}}) \}.$$
Denote by $M_{\mathfrak{p}}$ the quotient of $ P \cap
\C_{\mathfrak{p}}$ by the $S^1$-action. 
Then one can check the following facts: ${\mathcal{F}} = \{ M_{\mathfrak{p}}; \;
{\mathfrak{p}} \in {\mathcal{P}} \}$, the multiplicity of
$M_{\mathfrak{p}}$ is ${\mathfrak{p}}$, its complex dimension is $
\# I({\mathfrak{p}}) -1$. Furthermore, we have
$$ S^1_{ M_{\mathfrak{p}}} =  \{ \zeta \in \C^* ;\;  \zeta^{\mathfrak{p}} =1 \text{ and } \forall i \notin I( {\mathfrak{p}}), \; \zeta^{\mathfrak{p}_i} \neq 1 \},$$
Note that the $S^1_{ M_{\mathfrak{p}}}$ are mutually disjoint. Theorem \ref{int2} follows by using that 
 $$ G = \cup_{ {\mathfrak{p}} \in {\mathcal{P}} } S^1_{
 M_{\mathfrak{p}}}$$
and  if $\zeta \in S^1_{ M_{\mathfrak{p}}}$,
then $ c(\zeta ,M_{\mathfrak{p}}  )= \zeta$ and the $b_i
 (\zeta,M_{\mathfrak{p}}  )$ are the $\zeta^{{\mathfrak{p}}_i}$ with
 $i \notin I({\mathfrak{p}} )$.

For example, if $({\mathfrak{p}}_1,{\mathfrak{p}}_2,{\mathfrak{p}}_3)
= (2,4,3)$, then ${\mathcal{P}} = \{ 1,2,4,3 \}$. There are four
supports: $M_1 =M_r$, $M_2$ which is 1-dimensional and $M_4$, $M_3$
which consist of one point.
The subsets of $S^1$ associated to these supports are 
$$ S^1_1= \{ 1 \}, \; S^1_2= \{ -1  \},  \;  S^1_4=  \{ i , -i \}
\text{ and }    S^1_3= \{ e^{i \pi /3}, \;  e^{i 2\pi /3} \}. \qed $$
\end{rem}


\section{The kernel of the Szeg{\"o} projector of $M_r$} \label{prep}

In this section, we prove that the Szeg{\"o} projector
$\Pi_{r,k}$ of the orbifold $M_r$ is a
Fourier integral operator, which is the content of theorem \ref{BSor}. In the manifold case this result is a consequence of a theorem of Boutet de Monvel and Sj{\"o}strand on the kernel of the Szeg{\"o} projector associated to the boundary of a strictly pseudoconvex domain \cite{BoSj}.
 We won't adapt the proof of Boutet de Monvel and Sj{\"o}strand. Instead we will deduce this result from the same result for the Szeg{\"o} projector of $M$, i.e. we will show that 
$$ \Pi_k \in {\mathcal{F}} \Rightarrow  \Pi_{r,k} \in {\mathcal{F}}_r. $$
First we define the algebra $\fl$, that we introduced in section
\ref{vasavoir} for a free torus action. Then following  an idea of Guillemin and
Sternberg  (cf. appendix of \cite{GuSt1}), we prove that $\laP \in \fl$. Then we relate the algebras $\fl$ and $\fr$ as in
theorem \ref{keytheo2} and deduce theorem \ref{BSor}. 
To compare with
section \ref{vasavoir}, the proof of theorem \ref{FIO1} was 
in reverse order, i.e. we showed that $\Pi_{r,k} \in \fr
\Rightarrow \laP \in \fl$.

We think that the ansatz we propose for the kernel of
$\Pi_{r,k}$ is also valid for a general K{\"a}hler orbifold, not obtained
by reduction. But here, it is easier and more natural to deduce this
by reduction. 

\subsection{The algebra $\fl$}
Since the $\T^d$-action is not necessarily free, we have to modify
the definition of the algebra $\fl$. We consider $M$ as an orbifold and give local ansatz in orbifold charts as we did with the Toeplitz operators in section \ref{secork}. 

Introduce as in remark \ref{orbidata} an orbifold chart $(|U|,U,G,
\pi_U)$ of $M_r$ with the associated orbifold chart $(|V|,V,G,
\pi_V)$ of $P_\C$. Denote by $a_{g}$ (resp. $l_\te$) the action of $g\in G$
(resp. $\te \in \T^d$) on $V$. These actions both lift to $L^\ka_{V}
:= \pi_V^*L^\ka$. We denote them by ${\mathcal{A}}_g$ and ${\mathcal{L}}_\te$.

Define the data $(\Lambda_{V}, t^\ka_{\Lambda_V})$ corresponding to
the data $(\Lambda, t^\ka_{\Lambda})$ of section \ref{laFo}.  
\begin{gather*}
 \Lambda_{V} := \{ (\te,0,u,\te',0,u) \in V^2 ; \; \te,\te' \in \T^d \text{ and } u \in U \} \\
 t^{\kappa}_{\Lambda_V} (\te,0,u,\te',0,u)  := {\mathcal{L}}_{\te-
 \te'} .z \otimes z^{-1} \quad \text{if $z \in L_{V,(\te',0,u)}^{\ka}$
 and $z \neq 0$} 
\end{gather*}
Note that $t^{\kappa}_{\Lambda_V}$ is well-defined because the
$\T^d$-action on $V$ is free. 
In the definition of an operator $\op{T}_k$ of
$\fl$, we will assume that the lift $T_{k,VV}$ of its Schwartz kernel $T_k$ to
 $V^2$ satisfies
\begin{AS} \label{AsLa}
There exists a section $T'_{k,V}$ of $L^k_{V} \boxtimes L^{-k}_{V}$ invariant with respect to the action of $\T^d \times \T^d$ and the diagonal action of $G$ and 
of the form 
\begin{gather}  \label{AsLade}
 T'_{k,V}(x,y) = \Bigl( \frac{k}{2\pi} \Bigr)^{n-\frac{d}{2}} E^k_{\Lambda_V}(x,y) f(x,y,k) + O_\infty(k^{-\infty}) \end{gather}
on a neighborhood of $\Lambda_V$, where $E^\ka_{\Lambda_V}$ and
$f(.,k)$ satisfy the assumptions (\ref{kerFIO}.
{\romannumeral 1}) and (\ref{kerFIO}. {\romannumeral 2}) with $(\Ga, t_\Ga^\ka) = (\Lambda_{V}, t_{\Lambda_{V}}^\ka)$, such that  
\begin{gather*} 
 T_{k,VV} = \sum_{g \in G}  ({\mathcal{A}}_g \boxtimes \op{Id} )^* T'_{k,V} .\end{gather*}
\end{AS}

The whole definition of an operator of $\fl$ is the following. 
\begin{defin} 
$\fl$ is the set of operators $(\op{T}_k)$ with Schwartz kernel $T_k$ such that 
\begin{itemize}
\item $\Pi_{k}  \op{T}_k \Pi_{k} = \op{T}_k$ and ${\mathcal{L}}_\te^*
  \op{T}_k = \op{T}_k{\mathcal{L}}_\te^* = \op{T}_k$, for all $\te \in
  \T^d$.
\item  $T_k$ is  $O_{\infty}(k^{-\infty})$ on every
compact set $K \subset M^2$ such that $K \cap \Lambda = \emptyset$ 
\item $T_k$ satisfies assumption \ref{AsLa} for every orbifold chart $(|V|,V, G,\pi_V)$. 
\end{itemize}
\end{defin}

By adapting the proof of theorem \ref{symbeqFIO}, we define the
symbol map. 

\begin{prop} 
There exists a map $ \si : \fl \rightarrow \Ci (M_r)[[\hb]] $
which is onto and whose kernel consists of smoothing operators. 
\end{prop} 

The following theorem will be proved in the next subsection.
\begin{theo} \label{pilaFIO}
$\laP$ is an elliptic operator of ${\mathcal{F}}_\la$.
\end{theo}

\begin{rem} \label{comple}
 $\laP$ is elliptic means that its principal symbol doesn't
 vanish. This implies that $\laP$ itself doesn't vanish when $k$ is
 sufficiently large, so ${\mathcal{H}}_{\la,k}$ is not reduced to
 $(0)$ when $k$ is large enough. This completes the proof of the
 Guillemin-Sternberg theorem \ref{GuStTh} in the orbifold case. 
\qed \end{rem}

In proposition \ref{totsymb}, we defined the total symbol map 
$ \si : \fr \rightarrow \Ci(M_r)[[\hb]]$
and prove that its kernel consists of smoothing operators, without
using that $\Pi_{r,k} \in {\mathcal{F}}_r$. 
We can now generalize theorem \ref{keytheo2}.

\begin{theo} \label{keytheo3}
The map $\fl  \rightarrow \fr$ which sends $\op{T}_k$ into $\op{V}_k \op{T}_k \op{V}^*_k$ 
is  well-defined and bijective. Furthermore, the total symbols of
$\op{T}_k \in \fl$ and $\op{V}_k \op{T}_k \op{V}^*_k$ are the same. 
\end{theo}

The proof follows the same line as the proof of theorem
\ref{keytheo2}. To every chart  $(|V|,V,G,\pi_V)$ of $M$ is
associated a chart $(|U|,U,G,\pi_U)$. Assumption \ref{AsLa} corresponds to
assumption \ref{AsTo2}. 

As a corollary of this theorem, the total symbol map $\si : \fr
\rightarrow \Ci(M_r)[[\hb]]$ is onto. Furthermore it follows from the
stationary  phase lemma that ${\mathcal{F}}_r$ is a $*$-algebra and
the induced product $*_r$ on $\Ci(M_r)[[\hb]]$ is a star-product. 

Let $\tilde{\Pi}_{r,k}$ be an operator of ${\mathcal{F}}_r$ whose total symbol is the unit of $(\Ci(M_r)[[\hb]], *_r)$. 

\begin{lemme} \label{step2}
$ (\op{W}_k \tilde{\Pi}_{r,k} \op{W}_k^*)(\op{V}_k^*\op{V}_k) = \op{\Pi}_{\la,k} + \op{R}_k$
where $\op{R}_k$ is $O(k^{-\infty})$. 
\end{lemme}
\begin{proof} 
Since $\op{W}_k^* \op{V}_k^* = \Pi_{r,k}$ and $\op{V}_k^* \op{W}_k^* =\Pi_{\la,k}$, we have 
\begin{align*} 
  \op{W}_k \tilde{\Pi}_{r,k} \op{W}_k^* \op{V}_k^* \op{V}_k & =  \op{W}_k \tilde{\Pi}_{r,k}\op{V}_k \\
& =  \op{W}_k  \tilde{\Pi}_{r,k} \op{V}_k\op{V}_k^*   \op{W}_k^*
\end{align*}
By theorem \ref{keytheo3}, $ \op{V}_k\op{V}_k^* =\op{V}_k \Pi_{\la,k} \op{V}_k^*$ belongs to $\fr$ since $\Pi_{\la,k}$ belongs to ${\mathcal{F}}_\la$. From the symbolic calculus of the operators of ${\mathcal{F}}_r$, it follows that 
$$ \tilde{\Pi}_{r,k} (\op{V}_k\op{V}_k^* ) \equiv  (\op{V}_k\op{V}_k^* ) $$
modulo an operator of $\fr$ whose total symbol vanishes. Applying again theorem \ref{keytheo3}, we obtain that 
$$  \op{W}_k \bigl( \tilde{\Pi}_{r,k} (\op{V}_k\op{V}_k^* ) \bigr) \op{W}_k^* \equiv \op{W}_k \op{V}_k\op{V}_k^* \op{W}_k^*$$
modulo an operator of $\fl$ whose total symbol vanishes. The left-hand
side is equal to  $\op{W}_k \tilde{\Pi}_{r,k}
\op{W}_k^* \op{V}_k^*\op{V}_k$ and the right-hand side to $\Pi_{\la,k}$. This proves (\ref{step2}). \end{proof} 

Denote by $( \op{V}_k^* \op{V}_k )^{-1} $ the inverse of $ \op{V}_k^*
\op{V}_k$ on ${\mathcal{H}}_{\la,k}$, that is  
\begin{gather*}   
\Pi_{\la,k} ( \op{V}_k^* \op{V}_k )^{-1} \Pi_{\la,k} =  ( \op{V}_k^*
\op{V}_k )^{-1} \text{ and }  \\  ( \op{V}_k^* \op{V}_k ) ( \op{V}_k^*
\op{V}_k )^{-1}= ( \op{V}_k^* \op{V}_k )^{-1} ( \op{V}_k^* \op{V}_k )=
\Pi_{\la,k}. 
\end{gather*} 
Lemma \ref{step2} implies that 
\begin{gather} \label{lkjh}
  ( \op{V}_k^* \op{V}_k )^{-1} =   \op{W}_k \tilde{\Pi}_{r,k}
\op{W}_k^* - \op{R}_k ( \op{V}_k^* \op{V}_k )^{-1}.\end{gather} 
By theorem \ref{keytheo3}, $\op{W}_k \tilde{\Pi}_{r,k}
\op{W}_k^*$ belongs to $\fl$. 
\begin{lemme}
$ \op{R}_k ( \op{V}_k^* \op{V}_k )^{-1}$ is $O(k^{-\infty})$. 
\end{lemme}
\begin{proof} 
Since $\op{W}_k \tilde{\Pi}_{r,k} \op{W}_k^*$ belongs to $\fl$, its
kernel is $O(k^{n-\frac{d}{2}})$. So $\op{W}_k \tilde{\Pi}_{r,k}
\op{W}_k^*$ is $O(k^{n-\frac{d}{2}})$. Using that $\op{R}_k$ is
$O(k^{-\infty})$, it follows from \eqref{lkjh} that $ ( \op{V}_k^*
\op{V}_k )^{-1}$ is $O(k^{n-\frac{d}{2}})$.
\end{proof}
We deduce from this that  $( \op{V}_k^* \op{V}_k )^{-1} \in
\fl$.
 We have 
$$ \Pi_{r,k}= \op{V}_k  ( \op{V}_k^* \op{V}_k )^{-1} \op{V}_k^* .$$
Consequently theorem \ref{keytheo3} implies that $ \Pi_{r,k}  \in {\mathcal{F}}_r$.

\subsection{The projector $\Pi_{\la,k}$} \label{pilapila}

This section is devoted to the proof of theorem \ref{pilaFIO}. We use
that $\Pi_k \in
\mathcal{F}$ together with the following consequence of $\eqref{eqes}$
\begin{gather} \label{moy}
 \Pi_{ \la,k} ( \underline{x}, x ) =  \int_{\T^d} \bigr(
 ({\mathcal{L}}_{\te} \boxtimes \op{Id})^* \Pi_k \bigl) (\underline{x}, x) \; \de_{\T^d}(\te). \end{gather}
The only difficult point is to prove that $\Pi_{ \la,k}$
satisfies assumption \ref{AsLa}.

Introduce an orbifold chart $(|V|,V,G,\pi_V)$ of $M$ as in the previous
section. Denote by $\Pi_{ \la,k,VV}$ and $\Pi_{k,VV}$ the lifts of
$\Pi_{ \la,k}$ and $ \Pi_{k}$ to $V^2$. Then 
\begin{gather*}
 \Pi_{ \la,k,VV} ( \underline{v}, v ) =  \int_{\T^d} \bigl( ({\mathcal{L}}_{\te} \boxtimes \op{Id})^* \Pi_{k,VV}\bigr) (\underline{v}, v) \; \de_{\T^d}(\te). \end{gather*}
We know that $\Pi_k \in \mathcal{F}$. So $\Pi_{k,V}$ satisfies
assumption \ref{AsTo2} for the orbifold $M$. Denote by $\Pi'_{k,V}$ an
associated kernel such that  
\begin{gather*} 
 \Pi_{k,VV} = \sum_{g \in G}  ({\mathcal{A}}_g \boxtimes \op{Id} )^* \Pi'_{k,V} .\end{gather*}
If we prove that $\Pi_{ \la,k,V}'$ given by
 \begin{gather} \label{moyV}
 \Pi_{ \la,k,V}' ( \underline{v}, v ) =  \int_{\T^d}
 \bigl( ({\mathcal{L}}_{\te} \boxtimes \op{Id})^* \Pi_{k,V}' \bigr) (\underline{v},
 v) \; \de_{\T^d}(\te). \end{gather}
satisfies (\ref{AsLade}), we are done. So the proof is locally
 reduced to the manifold case. 

Let us relate the section $E^{\kappa}_{\Delta_V}$ and
$E^{\kappa}_{\Lambda_V}$ appearing in (\ref{kerFIO2}) and (\ref{AsLade}).
Let $$s^\ka : V \rightarrow L^\ka_V$$ be a holomorphic
$\T^d_{\C}$-invariant section which doesn't vanish. Introduce the real
function $H$ such that  
 $(s^\ka,s^\ka)(v) = e^{- \ka H(v)}$ and the unitary section 
$  t^\ka = e^{\ka H/2} s^\ka $
Let us write
\begin{gather*} 
 E_{\Delta_V}^{\ka}(\underline{v},v ) = e^{- \ka \phi_{\Delta}(\underline{v},v)} t^{\ka} (\underline{v}) \otimes t^{\ka} (v), \quad
 E_{\Lambda_V}^{\ka}(\underline{v},v ) = e^{- \ka \phi_{\Lambda}(\underline{v},v)} t^{\ka} (\underline{v}) \otimes t^{\ka} (v) 
\end{gather*} 
So we have
$$ \Pi'_{k,V}(\underline{v},  v ) = (\tfrac{k}{2\pi})^n e^{-k \phi_{\Delta} (\underline{v}, v)} f( \underline{v}, v,k) t^k(\underline{v}) \otimes t^k(v) + O_\infty(k^{-\infty}).$$
Assumption  (\ref{kerFIO}.{\romannumeral 1}) 
determines only the Taylor expansion of $E_\Gamma^{\ka}$ along $\Gamma$. Hence, the
functions $\phi_{\Delta}$ and $\phi_{\Lambda}$ are unique modulo a
function which vanishes to any order along the associated Lagrangian
manifold. 

Recall that we introduced a function $\ph$ in proposition
\ref{pcontrolnorm}. Let $\tilde{\ph} (
\underline{v},v)$ be a function such that such that $\tilde{\ph}  ( v,v) = \ph (v)$ and
$$\bar{\underline{Z}}\tilde{\ph} \equiv Z.\tilde{\ph} \equiv 0$$
modulo $O(| \underline{v} - v |^{\infty})$ for every holomorphic
vector field $Z$ of $V$. Here $\bar{\underline{Z}}\tilde{\ph}$
(resp. $Z.\tilde{\ph}$) denote the derivative of $\tilde{\ph}$
with respect to the vector field $(\bar{Z},0)$ of $V^2$ (resp. $(0,
Z)$). We use the same notation in the following. 
\begin{lemme}  \label{explo}
We can choose the functions $\phi_{\Lambda}$ and $\phi_{\De}$ in such a way that 
$$  \phi_{\De}(\underline{v},v) =   \phi_{\Lambda}(\underline{v},v) -
\tilde{\ph} ( \underline{v},v) $$
\end{lemme} 

\begin{proof} 
Recall that $V = \T^d \times \td \times U \ni ( \te , t, u)=v$. By
proposition \ref{pcontrolnorm}, we have 
\begin{gather} \label{relH} H (v) = H_r(u) + \ph (t,u) \end{gather}  
By reduction, $U$ is endowed with a complex structure
(cf. section \ref{comcomp}). Introduce a function $\tilde{H}_r(
\underline{u},u)$ such that $\tilde{H}_r(u,u) = H_r(u)$ and
$$\bar{\underline{Z}}\tilde{H}_r \equiv Z.\tilde{H}_r \equiv 0 \mod
O(| \underline{u} - u |^{\infty})$$ for every holomorphic vector field
$Z$ of $U$. Then 
\begin{gather} \label{poiu}
 \phi_{\Lambda}(\underline{v},v) := \tfrac{1}{2} (H(\underline{v}) +
H(v) ) - \tilde{H}_r(\underline{u},u) \end{gather} 
is a function associated to $\Lambda_V$. This is easily checked
using that 
$$  \nabla t^\ka = \frac{\ka}{2} ( \bar{\partial} H - \partial H)
\otimes t^\ka $$
Set $\tilde{H}(\underline{v},v) = \tilde{H}_r(\underline{u},u) +
\tilde{\ph} ( \underline{v},v)$. 
In the same way we get that 
$$ \phi_{\De}(\underline{v},v) :=  \tfrac{1}{2} (H(\underline{v}) +
H(v) ) - \tilde{H}(\underline{v},v) $$
is a function associated to $\De_V$. 
\end{proof}

\begin{lemme} We have over $V^2$
$$ \Pi'_{\lambda, k,V} (\underline{v},  v ) = \Bigl(
\frac{k}{2\pi}\Bigr)^{n- \frac{d}{2}}   e^{-k \phi_{\Lambda}(\underline{v},v)} g( \underline{v},  v,k)  t^k(\underline{v}) \otimes t^k(v) + O_\infty ( k^{-\infty})$$ 
where 
$$  g( \underline{v},  v,k) = \Bigl(\frac{k}{2\pi}\Bigr)^{\frac{d}{2}}\int_{W} e^{k \phi( \te', \underline{v},
  v )} f(\te+ \te',\underline{t}, \underline{u}, v,k) |d\te'| $$
$W$ is any neighborhood of $0$ in $\T^d$ and $\phi( \te', \underline{v}, v ) = \tilde{\ph}(\te+ \te',\underline{t}, \underline{u}, v)$.
\end{lemme}

\begin{proof} 
Since $ {\mathcal{L}}_{\te}^* t^k = t^k$,  (\ref{moyV}) implies 
$$ \Pi'_{\lambda, k,V} (\underline{v},  v ) \equiv \Bigl(
\frac{k}{2\pi} \Bigr)^n
t^k(\underline{v}) \otimes t^k(v) \int_{\T^d} e^{-k \phi_\De ( \underline{\te}+ \te',\underline{t}, \underline{u}, v)} f(\underline{\te}+ \te',\underline{t}, \underline{u}   , v,k) |d\te'| $$
modulo $O_\infty(k^{-\infty})$. 
Replacing $\te' + \underline{\te} - \te $ by $\te'$, this leads to 
$$ \Pi_{\lambda, k,V}' (\underline{v},  v ) \equiv
 \Bigl( \frac{k}{2\pi} \Bigr)^n
 t^k(\underline{v}) \otimes t^k(v) \int_{\T^d} e^{-k \phi_\De (\te+
   \te',\underline{t}, \underline{u}, v)} f(\te+ \te',\underline{t},
 \underline{u}   , v,k) |d\te'| $$
modulo $O_\infty(k^{-\infty})$.
Since the imaginary part of $\Phi_\De$ is positive outside the diagonal of $V^2$, we can restrict the integral over any neighborhood $W$ of $0$ in $\T^d$.
Now using lemma \ref{explo} and the fact that $\phi_\Lambda$ is
 independent of $\te$ which appears in equation \eqref{poiu}, we obtain the result. 
\end{proof}

To complete the proof of theorem \ref{pilaFIO}, it suffices to prove that $g (.,k)$ admits an asymptotic expansion in power of $k$ for the $\Ci$ topology on a neighborhood $\Lambda \cap V^2$. We prove this by applying the stationary phase lemma \cite{Ho1}. So the result is a consequence of the following lemma.  
\begin{lemme} \label{hyp}
Let $( \underline{v}_0 , v_0) \in \Lambda \cap V^2$. Then the Hessian $d^2_{\te'} \phi$ at $(0,\underline{v}_0,v_0)$ is a real definite positive matrix. Furthermore, on a neighborhood of $(0,\underline{v}_0,v_0)$ in $\T^d \times V^2$, 
$$ \phi = \textstyle{\sum} h_{ij}( \partial_{\te^{'i}} \phi ) ( \partial_{\te^{'j}} \phi ) $$
where the $h_{ij}$ are $\Ci$ functions of $\te',\underline{v},v$.
\end{lemme}

Before we prove this lemma, let us state some intermediate results. If $h$ is a function of $\Ci(V)$, we denote by $\tilde{h}$ a function of $\Ci(V^2)$ such that $\tilde{h}(\underline{v},v) = h(v)$ and 
\begin{gather} \label{rere}
\bar{\underline{Z}}. \tilde{h} \equiv Z.\tilde{h} \equiv 0 \mod O(| \underline{v} - v |^{\infty}) \end{gather}
for every holomorphic vector field $Z$ of $V$. 

 Denote by $t^i$ the coordinates of $t = \sum t^i \xi_i \in \td$. Let
 us compute the derivatives of $\tilde{h}$ with respect to the vector
 fields $\partial_{\te^j} = \xi^\#_j$ and $\partial_{t^j} = J\xi^\#_j$
 acting on the left and the right respectively.
\begin{lemme} \label{lextinvt}
If $h$ is $\T^d$-invariant, then $$i \partial _{\underline{\te}^j}\tilde{h} \equiv \partial_{\underline{t}^j}\tilde{h} \equiv -i \partial _{\te^j}\tilde{h} \equiv \partial_{t^j}\tilde{h} \equiv \tfrac{1}{2} \tilde{h}_j \mod O(| \underline{v} - v |^{\infty})$$
where $h_j= \partial_{t^j} h$.
\end{lemme}

In particular, we obtain the following relations 
\begin{gather} \label{relat} 
2\partial_{\underline{t}^j} \tilde{t}^i \equiv \de_{ij}, \quad 2 \partial_{t^j} \tilde{t}^i \equiv \de_{ij}, \quad 2 \partial_{\underline{\te}^j} \tilde{t}^i \equiv- i \de_{ij}, \quad 2\partial_{\te^j} \tilde{t}^i \equiv i \de_{ij}
\end{gather} 
modulo $ O(| \underline{v} - v |^{\infty})$.

\begin{proof} 
Since $[Z, \partial_{\te^j}]$ (resp. $[Z, \partial_{t^j}]$) is a holomorphic vector field when $Z$ is, the various derivatives we need to compute satisfy equations (\ref{rere}). So we just have to compute their restriction to the diagonal of $V^2$. To do this observe that  
$$  \langle d \tilde{h}, \partial_{\underline{\te}^j} + i \partial_{ \underline{t}^j} \rangle \eval{(v,v)} = 0, \qquad  \langle d \tilde{h}, \partial_{\te^j} - i \partial_{t^j} \rangle \eval{(v,v)} = 0$$
since  $\partial_{\te^j} - i \partial_{t^j} =  \xi^\#_j - iJ \xi^\#_j$ is holomorphic. Furthermore 
$$ \langle d \tilde{h}, \partial_{\underline{\te}^j} +  \partial_{ \underline{\te}^j} \rangle \eval{(v,v)} = 0, \qquad  \langle d \tilde{h}, \partial_{\underline{t}^j} +   \partial_{t^j} \rangle \eval{(v,v)} = \partial_{t^j}.h (v) $$
since $\tilde{h}(v,v) =h(v)$ and $h$ is $\T^d$-invariant. 
\end{proof}

\begin{proof}[Proof of lemma \ref{hyp}]  Let us write $ \ph(v) = \tfrac{1}{2} \textstyle{\sum} t^i t^j \ph_{ij}(v)$ on a neighborhood of $v_0$. Consequently, we have on a neighborhood of $(v_0,v_0)$
\begin{gather} \label{poi} 
\tilde{\ph}  =  \tfrac{1}{2} \textstyle{\sum} \tilde{t}^i \tilde{t}^j \tilde{\ph}_{ij}.
\end{gather}
It follows from equations (\ref{relat}) and proposition \ref{derivph} that 
\begin{align*} 
 \partial_{\underline{\te}^i} \partial_{\underline{\te}^j} \tilde{\ph} (v_0,v_0) =&  -\ph_{ij}(v_0) \\  =& -2 g( \xi^\#_i,\xi^\#_j)(v_0). \end{align*} 
The first part of the lemma follows. By (\ref{poi}) there exists $\Ci$ functions $h^1_{ij}$ such that 
$   \partial_{\underline{\te}^i} \tilde{\ph} = \textstyle{\sum} \tilde{t}^j h^1_{ij}$
on a neighborhood of $(v_0,v_0)$.
Derivating with respect to $\underline{\te}$, we obtain
$$h^1_{ij} (v_0,v_0) = -i \ph_{ij}(v_0).$$ Hence it is an invertible matrix and there exists functions $h^2_{ij}$ such that 
$    \tilde{t}^i = \textstyle{\sum} h^2_{ij}\partial_{\underline{\te}^j} \tilde{\ph} $ on a neighborhood of $(v_0,v_0)$. By (\ref{poi}),
$$  \tilde{\ph}   =  \textstyle{\sum} h^3_{ij} (\partial_{\underline{\te}^i} \tilde{\ph})( 
\partial_{\underline{\te}^j}\tilde{\ph} ) $$
for some $\Ci$ functions $h^3_{ij}$ on a neighborhood of $(v_0,v_0)$. The second part of the lemma follows. 
\end{proof}

\bibliography{biblio}

\end{document}